\documentclass[a4paper,11pt]{article}
\usepackage[noblocks]{authblk}
\usepackage{todonotes}
\usepackage{float}
\usepackage{amsmath}
\usepackage{amssymb}
\usepackage{amsthm}
\usepackage{pifont}
\usepackage{graphicx}
\usepackage{xcolor}
\usepackage{fullpage}
\usepackage[percent]{overpic}

\graphicspath{{./Figures/}}

\usepackage{xcolor}
\definecolor{newcolor}{rgb}{.8,.349,.1}
\usepackage{mathrsfs} 

\usepackage{amsmath}
\usepackage[noblocks]{authblk}
\usepackage{amsthm}
\usepackage{algpseudocode}
\usepackage{algorithm}
\usepackage{amsfonts}
\usepackage{graphicx}
\usepackage{geometry}
\usepackage{bm}
\usepackage{todonotes}
\usepackage{mathtools}
\usepackage{comment}
\usepackage{epstopdf}
\usepackage{float}
\usepackage{placeins}

\usepackage{xcolor}
\usepackage{amsmath}
\usepackage{amssymb}
\usepackage{amsthm}
\usepackage{amsfonts}
\usepackage{graphicx}
\usepackage{geometry}
\usepackage{bm}
\usepackage{todonotes}

\usepackage{mathtools}
\newtheorem{pro}{Proposition}

\usepackage[percent]{overpic}
\usepackage[latin1]{inputenc}
\usepackage{epstopdf}
\usepackage{float}
\usepackage{placeins}
\graphicspath{ {./Figures/} }
\newcommand{\td}[2]{\frac{ {\rm d} #1}{ {\rm d} #2}}
\newcommand{\BB}{\text{\Large$\mathbb{B}^\varepsilon$} }

\usepackage{booktabs} 
\usepackage{type1cm}        
\usepackage[noblocks]{authblk}
\usepackage{todonotes}
\usepackage{float}
\usepackage{amsmath}
\usepackage{amssymb}
\usepackage{amsthm}
\usepackage{pifont}
\usepackage{graphicx}
\usepackage{xcolor}                            
\usepackage{makeidx}         
\usepackage{graphicx}        
\usepackage{multicol}        
\usepackage[bottom]{footmisc}
\usepackage{booktabs}
\usepackage{siunitx}
\usepackage{newtxtext}       %
\usepackage[varvw]{newtxmath}       
\graphicspath{ {./Figures/} }
\newcommand{\pad}[2]{\frac{\partial #1}{\partial #2}}

\usepackage[percent]{overpic}

\title{Standard versus Asymptotic Preserving Time Discretizations for the Poisson-Nernst-Planck System in the Quasi-Neutral Limit}

\author{Clarissa Astuto}
\affil{Department of Mathematics and Computational Science, University of Catania, Italy}

\begin{document}
\maketitle

\begin{abstract}
In this paper, we investigate the correlated diffusion of two ion species governed by a Poisson-Nernst-Planck (PNP) system. Here we further validate the numerical scheme recently proposed in \cite{astuto2025asymptotic}, where a time discretization method was shown to be Asymptotic-Preserving (AP) with respect to the Debye length. 
For vanishingly Debye lengths, the so called Quasi-Neutral limit can be adopted, reducing the system to a single diffusion equation with an effective diffusion coefficient \cite{CiCP-31-707}. Choosing small, but not negligible, Debye lengths, standard numerical methods suffer from severe stability restrictions and difficulties in handling initial conditions. IMEX schemes, on the other hand, are proved to be asymptotically stable for all Debye lengths, and do not require any assumption on the initial conditions. In this work, we compare different time discretizations to show their asymptotic behaviors.
\end{abstract}

\section{Introduction}
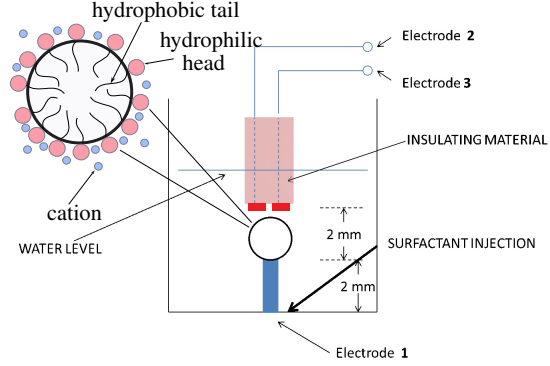
\begin{figure}[h]
	\centering
	\centering
\begin{overpic}[abs,width=0.6\textwidth,unit=1mm,scale=.25]{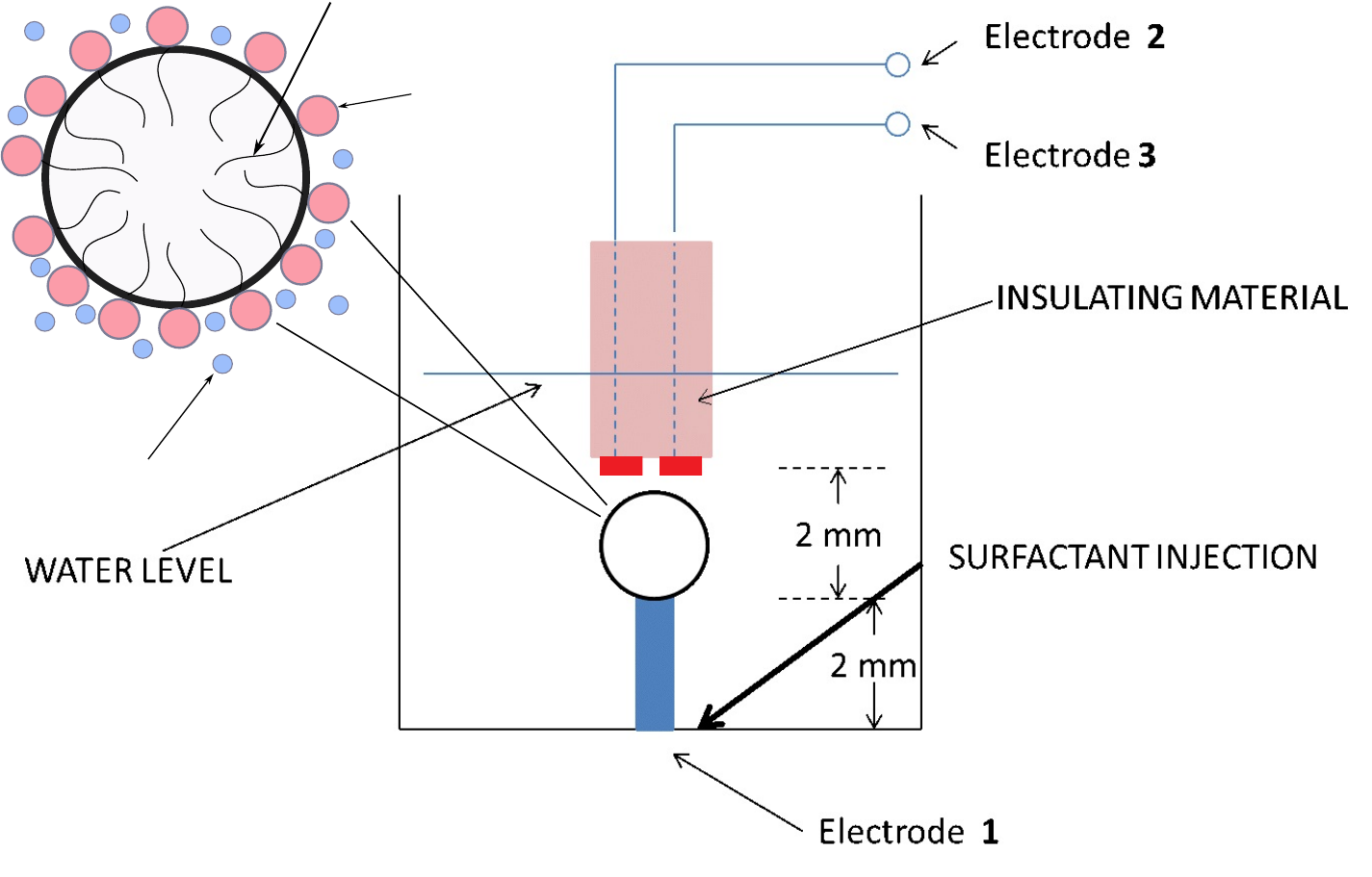}
\put(11,70){\footnotesize hydrophobic tail}
\put(30,63){\footnotesize hydrophilic}
\put(35,60){\footnotesize head}
\put(5,20){\footnotesize cation}
\end{overpic}
	\caption{\textit{Scheme of the experimental setup. The top-left inset shows the behavior of anions and cations at the bubble air-water interface: cations (blue) have hydrophilic heads, while anions (red), with their hydrophobic tails, result to be inside the bubble and their hydrophilic heads at the surface.}}
\label{fig_setup}
\end{figure}

Poisson-Nernst-Planck (PNP) system is a useful representation in many applications. The diffusion and migration of electric charge play a central role in a wide range of technologies and sciences \cite{lundstrom2002fundamentals,mason1988transport}: semiconductor technology controls the migration and diffusion of quasi-particles of charge in circuits \cite{markowich2012semiconductor,henisch1984semiconductor,nastasi2020full,nastasi2021efficient}; sheath formation refers to the development of a protective or structured layer around an object, material, or biological entity. The specific meaning depends on the context of plasma physics, when a boundary layer forms near a solid surface in contact with a plasma \cite{chen2012introduction,henderson1913fitness}. In this paper, we study the correlated motion of positive and negative ions near a (possible oscillating) bubble surface, modeling the cell membrane. The system consists of an anchored gas drop exposed to a diffusive flow of charged surfactants, which can reversibly adsorb onto the bubble surface (see Fig.~\ref{fig_setup}); \cite{Raudino20168574}. Their correlated diffusion is described by a Poisson-Nernst-Planck system, where the drift is driven by Coulomb interactions, obtained from a self-consistent Poisson equation.

An important physical parameter in this context is the Debye length, ($\lambda_D$, proportional to a small parameter) $\varepsilon$, that characterizes the distance over which electric fields are screened in a plasma due to local charge separation. While on large scales the plasma remains nearly quasineutral, with this short length scale charge imbalances can persist. Because the Debye length is often much smaller than the overall dimensions of the system, standard discretization strategies struggle to resolve it, and explicit time-stepping approaches become computationally very demanding \cite{degond2013asymptotic}.

Moreover, standard numerical schemes about the quasineutral limit assume that the initial conditions are well-prepared, meaning that they are compatible with the asymptotic model. In practice, however, this assumption is very restrictive for small values of $\varepsilon$  \cite{crouseilles2024high}.
In this paper, we design a numerical scheme that is independent of well prepared initial conditions in the limit as the Debye length approaches zero.

In a fluid, the time evolution of the concentrations of cations $c_+= c_+(\vec{x},t)$ and of anions $c_-= c_-(\vec{x},t)$, in presence of a electrostatic potential is governed by the following conservation laws
\begin{subequations}
\label{eq_full_model}
\begin{align}
		\displaystyle \frac{\partial c_\pm}{\partial t} &= D_\pm \Delta c_\pm \pm \chi_\pm \nabla \cdot \left(c_\pm \nabla ( q\varphi)\right), \quad {\rm in }\, \Omega \\
		\displaystyle J_\pm &= - D_\pm \nabla c_\pm \mp \chi_\pm c_\pm \nabla ( q\varphi), \quad {\rm in }\, \Omega\\ 
		-\bar \varepsilon \Delta (q\varphi) &= \frac{c_+}{\tilde m_+} - \frac{c_-}{\tilde m_-}, \quad {\rm in }\, \Omega\\		
		\displaystyle \nabla c_\pm \cdot \widehat n &= 0, \quad {\rm on }\, \partial\Omega \\
        \nabla \varphi \cdot \widehat n &= 0, \quad {\rm on }\, \partial\Omega.
\end{align}
\end{subequations}
where $J_\pm$ are the fluxes, $D_\pm$ are the diffusion coefficients respectively for $c_\pm$, $\chi_\pm = D_\pm/k_BT$ are the ions mobilities, where $k_B$ is the Boltzmann's constant and $T$ is the absolute temperature (assumed to be constant). $\varphi$ is the electrostatic potential that accounts for the interaction among the diffusing ions, $q $ is the (positive) electron charge, $\tilde m_\pm $ are the molecular mass of positive and negative ions, respectively, and $\bar \varepsilon$ is a quantity that is proportional to the Debye length. 

Following the dimensional analysis we performed in \cite{CiCP-31-707}[Appendix.2], Eqs.~\eqref{eq_full_model} can be rewritten as follows
\begin{subequations}
\label{eq_full_model_adim}
\begin{align} \label{eq_full_model_adim_cpm}
		\displaystyle \frac{\partial c_\pm}{\partial t} &= D_\pm \left(\Delta c_\pm \pm \nabla \cdot \left(c_\pm \nabla \Phi\right)\right), \quad {\rm in }\, \Omega\\ \label{eq_full_model_adim_Phi}
		-\varepsilon \Delta \Phi &= \frac{c_+}{m_+} - \frac{c_-}{m_-}, \quad {\rm in }\, \Omega\\		
		\displaystyle \nabla c_\pm \cdot \widehat n &= 0, \quad {\rm on }\, \partial\Omega \\
        \nabla \Phi \cdot \widehat n &= 0, \quad {\rm on }\, \partial\Omega,
\end{align}
\end{subequations}
where $\Phi = q\varphi$ and $\varepsilon = K^{-1}$, another quantity that is proportional to the Debye length, where 
\[ \displaystyle K = \frac{q^2N_A \rho}{\epsilon_0\epsilon_r m_0 k_B T},\] with $N_A $ the Avogadro's number, $\tilde{m}^\pm = m_0 m_\pm$ the molar mass of ions 
(Kg/mol) and $\rho^\pm$ their mass densities  {(Kg/m$^3$)}, with the assumption that $\rho:= \rho^+ = \rho^-$.

In this work, we adopt the same domain as in \cite{astuto2025asymptotic}. However, our primary focus is on the numerical challenges posed by the stiffness coming from the negligible values of the Debye length. In contrast to \cite{astuto2025asymptotic}, where a multiscale analysis was conducted to take into account the external potentials that describe the relation between the bubble surface and the ions, in this paper, we simplify the analysis by disregarding these additional scales. To this end, we impose homogeneous Neumann boundary conditions on $\partial \Omega$ (instead of the time-dependent boundary conditions defined in \cite{astuto2025asymptotic}), allowing us to isolate and better understand the effects of the Debye length on the system's behavior. 

The objective of this paper is to demonstrate that the IMEX strategy, when applied to the formulation described in Section~\ref{sec:QNL}, ensures the expected asymptotic properties, such as stability and accuracy. To this end, we compare different semi-implicit time-integration strategies, first applied to the original problem in~\eqref{eq_full_model_adim}, and later on to the reformulated system in~\eqref{eq_CQP_system}. In addition, we also investigate higher-order IMEX schemes in time, but only for academic purposes, since the spatial discretization remains second order.

\section{Quasi--Neutral Limit (QNL) }
\label{sec:QNL}In this section, we focus on the limit of the Debye length that goes to zero. It is well known that solving the coupled Poisson-Nernst-Planck (PNP) system \cite{EISENBERG2007,LU20112475,CiCP-31-707} gives a computational cost that is prohibitively high due to the stability requirements imposed by the Debye length. Our goal is to reformulate the system~\eqref{eq_full_model_adim} in a way that allows the construction of a numerical scheme free from time-step restrictions.

Expanding the derivatives in Eq.~\eqref{eq_full_model_adim_cpm}, we obtain
\begin{equation}
\label{eq:stiff_cpm}
    \pad{c_\pm}{t} = D_\pm\left(\Delta c_\pm \pm \nabla c_\pm \cdot \nabla \Phi \mp \frac{c_\pm}{\varepsilon} \left(\frac{c_+}{m_+} - \frac{c_-}{m_-} \right) \right)
\end{equation}
where we made use of $\displaystyle -\Delta \Phi = \frac{1}{\varepsilon}\left(\frac{c_+}{m_+} - \frac{c_-}{m_-}\right)$, a term that generates an indeterminate form in the equations when $\varepsilon \to 0$. Classical discretization methods fail to capture the small scales corresponding to the Debye length, which can be much smaller than the overall size of the system. Using explicit time discretization schemes in such cases leads to prohibitively expensive simulations, due to the restriction in the time step \cite{degond2013asymptotic,CiCP-31-707}. For this reason, we propose a new formulation of problem \eqref{eq_full_model_adim} that follows the strategy used to obtain the so called \textit{Quasi--Neutral Limit} (QNL) regime (see, for instance, \cite{jungel,CiCP-31-707}).

We rewrite system~\eqref{eq_full_model_adim}, copying the equation for $\Phi$, and defining two new quantities 
\[ \mathcal Q = \frac 1 \varepsilon\left(\frac{c_+}{m_+} - \frac{c_-}{m_-}\right), \qquad \mathcal C = \frac{c_+}{m_+} + \frac{c_-}{m_-},\]
obtaining
\begin{subequations}
\label{eq_CQP_system}
\begin{align}
\label{eq_CQP_system_C}
    \frac{\partial \mathcal C}{\partial t} & = \widetilde D \Delta \mathcal C + \varepsilon \widehat D \Delta \mathcal Q + \nabla \cdot \left(\left(\widehat D \mathcal C + \widetilde D \varepsilon \mathcal Q \right)\nabla \Phi \right) \qquad {\rm in } \, \Omega \\
    \label{eq_CQP_system_Q}
  \frac{\partial \mathcal Q}{\partial t} & =  \frac{\widehat D}{\varepsilon} \Delta \mathcal C + \widetilde D \Delta \mathcal Q+ \nabla \cdot \left(\left( \frac{\widetilde D}{\varepsilon} \mathcal C + \widehat D \mathcal Q \right)\nabla \Phi \right) \qquad {\rm in } \, \Omega \\ \label{eq_CQP_system_Phi}
   -\Delta \Phi & = \mathcal Q \qquad {\rm in } \, \Omega\\
  \nabla \mathcal C \cdot \hat n & = \nabla \mathcal Q \cdot \hat n = \nabla \Phi \cdot \hat n = 0 \qquad {\rm on } \, \Gamma 
\end{align}
\end{subequations}
where $\widetilde D = \left({D_+ + D_-}\right)/{2}$ and $ \widehat D = \left({D_+ - D_-}\right)/{2}$ are the new diffusion coefficients. 
While the goal of~\cite{astuto2025asymptotic} was to design a time discretization of system~\eqref{eq_CQP_system} that is uniformly stable with respect to $\varepsilon \to 0$ and consistent with the time discretization of the limit $\varepsilon = 0$, here we propose different time discretizations and compare their performance in order to identify which scheme is more suitable.

\section{Numerical schemes} 
In this section, we describe the space and time discretization. We follow the strategy introduced in \cite{astuto2025asymptotic} where a second order fully discrete scheme for the Multiscale Poisson-Nernst-Planck system was developed. For this reason, we defer the detailed description to the Appendix. The main novelty is the comparison of IMEX schemes applied to two different formulations of the same system. In both cases, we solve the linear systems in a monolithic way. Secondly, we compare these techniques to a standard one that solves the unknowns in different steps. The extension to higher-order IMEX schemes is mainly of academic interest, as the spatial discretization remains second order.

\subsection{Space discretization and variational formulation}
\label{sec:space_discretization}
Here we describe the spatial discretization. The two time-dependent equations include both diffusion and drift terms. When the mesh P\'{e}clet number (pec) is high, standard numerical methods may lead to instabilities. In general, for finite difference schemes, Wesseling \cite{Wesseling2023600} shows that numerical solutions remain stable if the mesh P\'{e}clet number is kept below two in each spatial dimension. The issue of spatial oscillations has also been extensively studied in the context of finite element methods \cite{christie1976finite,heinrich1977upwind,heinrich1977quadratic,smith1980finite,brezzi1992relationship}. In \cite{smith1980finite}, the author studies the case for finite element methods using both linear and quadratic symmetric basis functions, and forms a comparison of the various numerical schemes. It is shown that the numerical solutions produce spatial oscillations whenever the mesh P\'{e}clet number exceeds two.

We make use of a ghost nodal finite elements method, a recently developed numerical scheme that achieves second order  accuracy in space \cite{astuto2025nodal,astuto2024comparison}. The same ghost-FEM discretization has been extended to the numerical solution of biological network formation in a leaf-shaped domain in \cite{astuto2024self}. In \cite{dilip2025multigrid}, the authors present multigrid methods for solving elliptic partial differential equations on arbitrary domains using the same ghost-FEM discretization.

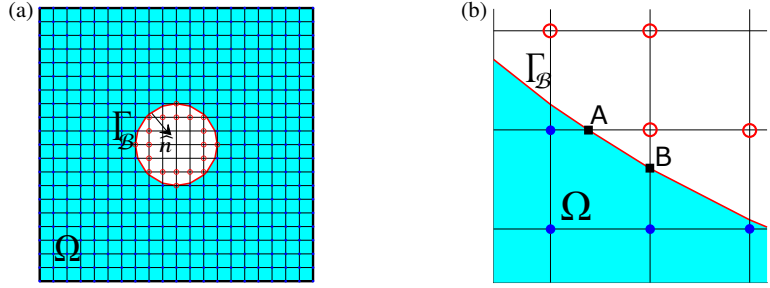
\begin{figure}[h]
    \centering
    \begin{minipage}{.49\textwidth}
\centering\begin{overpic}[abs,width=0.75\textwidth,unit=1mm,scale=.25]{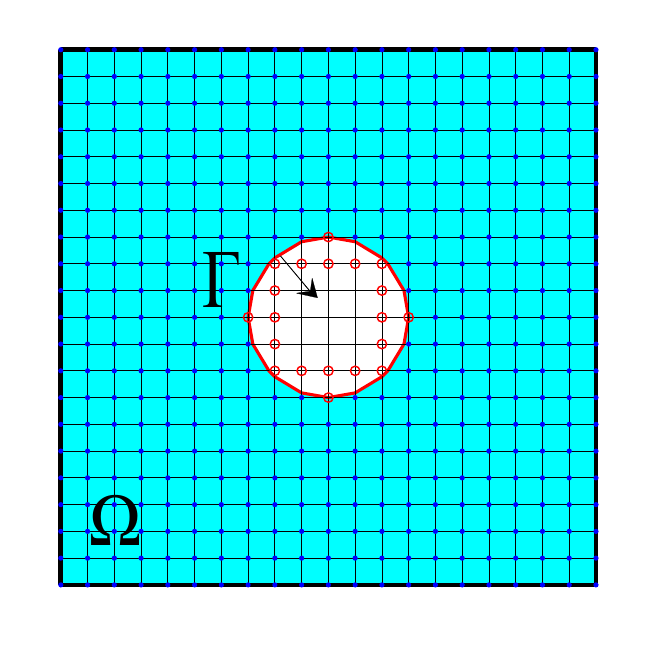}        
\put(0,40){(a)}
\put(20,22){$\widehat n$}
\put(14.3,22.5){$\mathcal B$}
\end{overpic}
    \end{minipage}
    \begin{minipage}{.49\textwidth}
\centering\begin{overpic}[abs,width=0.75\textwidth,unit=1mm,scale=.25]{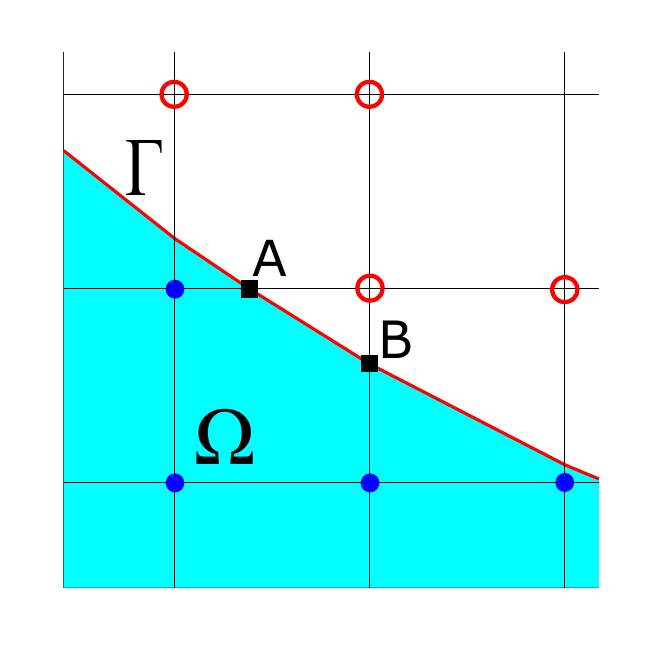}        
\put(0,40){(b)}
\put(9.6,30.5){$\mathcal B$}
\end{overpic}
    \end{minipage}
\caption{\textit{Discretization of the computational domain. $\Omega$ is the light blue region inside the unit square $R$. (a): classification of the grid points: the blue points are the internal ones while the red circles denote the ghost points. (b): points of intersection between the grid and the circular boundary $\Gamma_\mathcal{B}$.}}  
\label{fig:Domain2D}
\end{figure}

In this paper, the domain under consideration is a square with a circular hole in the middle. We start considering the square region $R = [-L_x/2,L_x/2]\times [-L_y/2,L_y/2]$ and the domain $\Omega = R\setminus \mathcal{B}$, where $\mathcal{B}$ is the circle centered at the origin, with radius $r_{\mathcal{B}}$ (see Fig.~\ref{fig:Domain2D} (a)). The domain is the one introduced in \cite{astuto2025asymptotic}, but here we specifically focus on the stiffness arising from the Debye length, while avoiding the multiple scales associated with time-dependent boundary conditions previously considered. Therefore, we assume homogeneous Neumann boundary conditions on $\Gamma := \partial \Omega$.

Following the approach showed in  \cite{Osher2002,Russo2000,book:72748, Sussman1994}, the domain $\Omega$ is implicitly defined by a level set function $\phi(x,y)$ that is negative inside $\Omega$, positive in $R\setminus \Omega$ and zero on the boundary $\Gamma_\mathcal{B} \subset \Gamma$, such that:
\begin{eqnarray}
	\Omega = \{(x,y): \phi(x,y) < 0\}, \qquad
	\Gamma_\mathcal{B} = \{(x,y): \phi(x,y) = 0\}.
\end{eqnarray}
We remark that we also impose zero Neumann boundary conditions in $\Gamma_\mathcal{B} \subset \Gamma.$
Now we introduce the following spaces of functions $V$ and $Q$
\begin{equation}
    V = \biggl\{ v\in H^1(\Omega) \biggr\}, \qquad Q = \biggl\{ q\in H^1(\Omega) : \int_\Omega q\,{ \rm d}\Omega = 0\biggr\}.
\end{equation}
As mentioned in the Introduction, to compare different semi-implicit discretizations, we examine two different formulations of our problem. Specifically, two different variational formulations, beginning with the one corresponding to the system presented in~\eqref{eq_full_model_adim}. 
\begin{pro}
Find $c_\pm \in V$ and $\Phi \in Q$ such that  for almost every $t\in(0,T)$, it holds
\begin{subequations}
\label{pro:variational1}
\begin{align}
\left( \frac{\partial c_+}{\partial t}, v\right)_{L^2(\Omega)} = & - D_+ \left(\left( \nabla c_+,\nabla v \right)_{L^2(\Omega)} + \left(c_+ \nabla \Phi,\nabla v \right)_{L^2(\Omega)} \right) \\ 
\left( \frac{\partial c_-}{\partial t}, u\right)_{L^2(\Omega)} = & - D_- \left(\left( \nabla c_-,\nabla u \right)_{L^2(\Omega)} - \left(c_- \nabla \Phi,\nabla u \right)_{L^2(\Omega)} \right) \\
\varepsilon \left( \nabla \Phi,\nabla q\right)_{L^2(\Omega)} = &\left(\frac{c_+}{m_+} - \frac{c_-}{m_-},q\right)_{L^2(\Omega)}
\end{align}
\end{subequations}
\end{pro}
Secondly, we show the variational formulation of the system showed in~\eqref{eq_CQP_system}.
\begin{pro}
Find $\mathcal C, \mathcal Q \in V$ and $\Phi \in Q$ such that for almost every $t\in(0,T)$, it holds
\begin{subequations}
\label{pro:variational}
\begin{align}
\left( \frac{\partial \mathcal{C}}{\partial t}, v\right)_{L^2(\Omega)} = & - \widetilde D \left( \nabla \mathcal{C},\nabla v \right)_{L^2(\Omega)}  - \varepsilon \widehat D \left( \nabla \mathcal Q, \nabla v \right)_{L^2(\Omega)} \\ \nonumber & - \left(\left( \widehat D \mathcal C + \varepsilon \widetilde D \mathcal Q \right) \nabla \Phi,\nabla v \right)_{L^2(\Omega)}  \\  
{\varepsilon}\left( \frac{\partial \mathcal{Q} }{\partial t}, u \right)_{L^2(\Omega)} = & - {\widehat D} \left( \nabla \mathcal{C} ,\nabla u  \right)_{L^2(\Omega)}  -  {\varepsilon}\widetilde D \left( \nabla \mathcal Q , \nabla u  \right)_{L^2(\Omega)} \\ \nonumber& - \left(\left( {\widetilde D} \mathcal C  +  {\varepsilon}\widehat D \mathcal Q  \right) \nabla \Phi ,\nabla u  \right)_{L^2(\Omega)}  \\ 
\left( \nabla \Phi,\nabla q\right)_{L^2(\Omega)} = &\left(\mathcal{Q},q\right)_{L^2(\Omega)}
\end{align}
\end{subequations}
\end{pro}

The discretization of the problem is obtained by discretizing the rectangular region $R\subset \mathbb R^2$ (with $\Omega \subset R$) by a regular rectangular grid (see Figure~\ref{fig:Domain2D}).
{Now we introduce the family of finite-elements spaces 
\begin{subequations}
\label{eq:basis}
\begin{align}\label{eq:Vh}
V_h & = \{ v_h \in V : {v_h}|_K \in \mathbb Q_1(K)\quad \forall K\in \Omega_h\} \\ \label{eq:Qh}
Q_h & = \{ q_h \in Q : {q_h}|_K \in \mathbb Q_1(K)\quad \forall K\in \Omega_h\}
\end{align}
\end{subequations}
where $\mathbb Q_1(K)$ denotes the space of piecewise bilinear functions on $K$,} 
which are continuous in $R$. To solve the variational problems (\ref{pro:variational1}-\ref{pro:variational}), we employ a finite-dimensional discretization. Specifically, the functions $c_\pm, \mathcal C, \mathcal Q \in V$ and $\Phi \in Q$ are  approximated by functions $c_{\pm,h}, \mathcal C_h, \mathcal Q_h$ and $\Phi_h$, that belong to finite-dimensional subspaces $V_h$ and $Q_h$, respectively.  
To perform our computations, the domain $\Omega$ is approximated by a polygonal domain $\Omega_h$. This approximation also extends to the boundary $\Gamma$, which is represented by $\Gamma_h$, {the straight line segment connecting the points of intersection between $\Gamma$ and the boundary of the tessellation}. Consequently, the original integrals defined over $\Omega$ and its boundary $\Gamma$ are now evaluated over $\Omega_h$ and $\Gamma_h$, respectively. A detailed description of the finite element discretization is provided in Appendix.

To reformulate the systems (\ref{pro:variational1}-\ref{pro:variational}) using the discrete spaces and the computational matrices for the spatial derivatives, we write
\begin{subequations}
\label{eq:FEM_semidiscr_1}
\begin{align}
\mathbb B[v_h]\frac{\partial c_{+,h}}{\partial t} & = - D_+ \left( \mathbb L[v_h] \, c_{+,h} + \mathbb H\left[ c_{+,h},v_h \right]  \Phi_h \right)  \\ 
 \mathbb B[u_h] \frac{\partial c_{-,h}}{\partial t} & = - D_- \left( \mathbb L[u_h]\, c_{-,h} - \mathbb H \left[c_{-,h},u_h \right]  \Phi_h \right) \\
\varepsilon\mathbb L[q_h]\, \Phi_h & = \mathbb B[q_h] \, \left( \frac{c_{+,h}}{m_+} - \frac{c_{-,h}}{m_-}\right)
\end{align}
\end{subequations}
and
\begin{subequations}
\label{eq:FEM_semidiscr}
\begin{align}
\mathbb B[v_h]\frac{\partial \mathcal{C}_h}{\partial t} & = - \widetilde D \, \mathbb L[v_h] \, \mathcal{C}_h  - \varepsilon \widehat D\, \mathbb L[v_h]\, \mathcal Q_h - \mathbb H\left[ \widehat D \mathcal C_h + \varepsilon \widetilde D \mathcal Q_h,v_h \right]  \Phi_h  \\ 
\varepsilon \mathbb B[u_h] \frac{\partial \mathcal{Q}_h}{\partial t} & = - {\widehat D}\, \mathbb L[u_h]\, \,\mathcal{C}_h  -  \varepsilon \widetilde D \,\mathbb L[u_h]\, \, Q_h - \mathbb H \left[ {\widetilde D}\mathcal C_h + \varepsilon \widehat D \mathcal Q_h,u_h \right]  \Phi_h \\
\mathbb L[q_h]\, \Phi_h & = \mathbb B[q_h] \, \mathcal{Q}_h 
\end{align}
\end{subequations}
where $\mathbb L[v_h] $ is the discrete operator that defines the stiffness matrix $\left( \nabla \diamond,\nabla v_h \right)_{L^2\left(\Omega_h\right)}$, such that $ \mathbb L[v_h] \, \mathcal C_h \approx \left( \nabla \mathbb \mathcal C_h,\nabla v_h\right)_{L^2\left(\Omega_h\right)}$; $ \mathbb H[\mathcal{C}_h,v_h]$ is the operator such that $\mathbb H \left[ \mathcal C_h,v_h \right]  \Phi_h \approx \left(  \mathcal C_h \nabla \Phi_h,\nabla v_h  \right)_{L^2\left(\Omega_h\right)}$; finally, we define the discrete operator $\mathbb B[v_h]$ for the mass matrix, such that $\mathbb B[v_h] \, \mathcal C_h \approx \left(\mathcal C_h,v_h\right)_{L^2\left(\Omega_h\right)}$.

Regarding the stability restriction due to the mesh P\'{e}clet number, we do the following considerations. From the definition in one dimension in \cite{Wesseling2023600}, we must satisfy
\begin{equation}
\label{eq:pec}
   \left| \partial_x \Phi  \right| < 2/h. 
\end{equation}
Considering now Eq.~\eqref{eq_full_model_adim_Phi} in one dimension as well, we have
\begin{equation} 
\label{eq:pec2}
    \left| \partial_x \Phi\right| \lesssim h\varepsilon^{-1} \frac{c_+}{m_+}
\end{equation}
where we considered the approximation $ \left| \partial_x \partial_x \Phi\right| \approx \frac 1 h  \left| \partial_x \Phi\right|.$ Combining Eqs.(\ref{eq:pec}-\ref{eq:pec2}), we have the following condition to satisfy at each time step
\begin{equation}
    \label{eq:pec3}
    \frac{c_\pm}{m_\pm} < 2 h^{-2} \varepsilon.
\end{equation}

\subsection{Time discretization}
We now discuss the discretization in time. The main objective is to design a semi-implicit time integrator that remains stable for all values of $\varepsilon$ (AP property) while preserving the order of accuracy in the asymptotic limit (AA property). The concept of an Asymptotic Preserving method has been introduced in \cite{jin1999efficient}. The commutative diagram below shows the main property of AP schemes, defined as follows. A scheme $\mathcal P^\varepsilon_h$ for $\mathcal{P}^\varepsilon$ with discretization parameters $h$ is called Asymptotic Preserving (AP) if it is stable independently from $\varepsilon$, and it converges to a numerical scheme $\mathcal P^0_h$ which is a numerical discretization of the problem $\mathcal P^0$ for a fixed discretization parameters $h$. 

\begin{figure}[h]
{\Large \[
  \begin{matrix}
\mathcal P_h^\varepsilon & {\longrightarrow} & \mathcal P^\varepsilon\\
\downarrow &  &  \downarrow\\
\mathcal P^0_h & \longrightarrow & \mathcal P^0 \\
\end{matrix}
\]
}

\caption{The AP diagram. $\mathcal P^\varepsilon$ is the original problem and $\mathcal P^\varepsilon_h$ its numerical approximation characterized by a discretization parameter $h$. The AP property corresponds to the request that $\mathcal P^0_h$ is consistent with $\mathcal{P}^0$ as $\varepsilon \to 0$, independently of $h$.}
\end{figure}

AP schemes are extremely powerful tools as they allow the use of the same scheme to discretize $\mathcal{P}^\varepsilon$ and $\mathcal{P}^0$, with fixed discretization parameters. Common numerical schemes do not provide this property, and tend to develop instabilities when the microscopic scale is under-resolved, for example, for negligible values of $\varepsilon$.


To introduce a time discretization, we consider a final time $T$ and define the time step as $\Delta t = T/N_{\rm ts},\, N_{\rm ts} \in \mathbb N,$ denoting the nodes in time by $t^n = n\Delta t$ and $\mathcal{C}_h^n \approx \mathcal{C}_h(t^n), \, n = 0,\cdots,N_{\rm ts}$. A semi--implicit discretization is adopted to achieve second and third order of accuracy in time. In particular, we make use of implicit-explicit (IMEX) Runge--Kutta schemes \cite{pareschi2000implicit,boscarino2016high,astuto2023self,boscarino2024implicit}, which are multi-step methods based on $s$-stages. 

We rewrite Eqs.~(\ref{eq:FEM_semidiscr_1}-\ref{eq:FEM_semidiscr}) in a vectorial form
\begin{equation}
\label{eq:vectorial}
 \BB\td{\bf Q}{t} =\Theta [{\bf Q}] {\bf Q}
\end{equation}
where we distinguish the two systems: system~\eqref{eq:FEM_semidiscr_1} becomes
\begin{equation}
\label{eq:vectorial1}
 \BB  = \begin{pmatrix}
 \mathbb B[v_h] & \underbar 0 & \underbar 0 \, \\
\underbar 0 & \mathbb B[u_h] &  \underbar 0 \,\\
\underbar 0 & \underbar 0 & \underbar 0 \, \\
\end{pmatrix} ,\qquad \Theta[{\bf Q}] = \begin{pmatrix}
- D_+ \mathbb L[v_h] & \underbar 0 & - D_+ \mathbb H\left[c_{+,h,}v_h\right]\\
\underbar 0 & - D_- \mathbb L[u_h] &  D_-\mathbb H\left[c_{-,h},u_h\right]\\
\mathbb B[q_h] & - \mathbb B[q_h] & \varepsilon \mathbb L[q_h] \\
\end{pmatrix} 
\end{equation}
where ${\bf Q} = [c_{+,h}\, ,\, c_{-,h}\, ,\, \Phi_h]^\top$. 
While system~\eqref{eq:FEM_semidiscr} becomes
\begin{equation}
 \BB  = \begin{pmatrix}
 \mathbb B[v_h] & \underbar 0 & \underbar 0 \, \\
\underbar 0 & \varepsilon \mathbb B[u_h] &  \underbar 0 \,\\
\underbar 0 & \underbar 0 & \underbar 0 \, \\
\end{pmatrix} ,\qquad \Theta[{\bf Q}] = \begin{pmatrix}
- \widetilde D \mathbb L[v_h] & - \varepsilon \widehat D \mathbb L[v_h] & - \mathbb H\left[\widehat D \mathcal C_h + \widetilde D \varepsilon \mathcal Q_h,v_h\right]\\
- {\widehat D} \mathbb L[u_h] & - \varepsilon\widetilde D \mathbb L[u_h] &  - \mathbb H\left[{\widetilde D} \mathcal C_h + \widehat D \varepsilon \mathcal Q_h,u_h\right]\\
\underbar 0 & - \mathbb B[q_h] & \mathbb L[q_h] \\
\end{pmatrix} 
\label{eq:vectorial2}
\end{equation}
and ${\bf Q} = [\mathcal C_h,\mathcal Q_h,\Phi_h]^\top$. 

Let us first set $\textbf{Q}^1_E = \textbf{Q}^n$, then the stage fluxes are calculated as
\begin{subequations}  
\label{eq_imex}
\begin{align}
\label{eq_imex_QE}
    \BB  \textbf{Q}_E^{i}& =  \BB \textbf{Q}^n + \Delta t\,\sum_{j=1}^{ i-1}\widetilde a_{i,j}\Theta(\textbf{Q}_E^j)\textbf{Q}_I^j, \quad i = 1,\cdots,s \\
\label{eq_imex_QI}
     \BB \textbf{Q}_I^{i} &=  \BB  \textbf{Q}^n + \Delta t\,\sum_{j=1}^{ i}  a_{i,j}\Theta(\textbf{Q}_E^j)\textbf{Q}_I^j, \quad i = 1,\cdots,s 
\end{align}
\end{subequations}
and the numerical solution is finally updated with
\begin{align}
\label{eq_imex_Q_bi}
    \BB  \textbf{Q}^{n+1} &=  \BB  \textbf{Q}^n + \Delta t\sum_{i=1}^{\rm s}  b(i) \Theta(\textbf{Q}_E^i)\textbf{Q}_I^i.
\end{align}
In our numerical tests, ${\rm s} = 2,3,4$. 

\subsubsection{Comparison of IMEX schemes: accuracy and stability analysis}
We employ several IMEX Runge-Kutta schemes to integrate the systems efficiently. These schemes differ in their stability properties (A-stable or L-stable) and in their ability to preserve the asymptotic limit as $\varepsilon \to 0$. An additive SA property allows the scheme to preserve its order of accuracy even in this limit.

\subsubsection*{I1: IMEX-H(2,2,2)}
We start with a second order, two stages SDIRK IMEX-RK methods of type I \cite{ascher1997implicit,pareschi2003high,pareschi2005implicit}, with the double Butcher tableau of the form
\begin{equation}
    \begin{array}{c|cc}
        0 & 0 & 0  \\
        1 & 1 & 0 \\ \hline
         & 1/2 & 1/2 
    \end{array} \hspace{4cm}
        \begin{array}{c|cc}
        1/2 & 1/2 & 0  \\
        1/2 & 0 & 1/2 \\ \hline
         & 1/2 & 1/2 
    \end{array}
\end{equation}
This scheme is the combination of the Heun method and a second-order A-stable DIRK
method. We refer to this scheme as IMEX-H(2,2,2).
\subsubsection*{I2: IMEX-SA(2,2,2)}
Second order numerical scheme. The implicit part is L-stable and stiffly accurate
\begin{equation}
    \begin{array}{c|cc}
        0 & 0 & 0  \\
        1 & 1/(2\gamma) & 0 \\ \hline
         & 1/2 & 1/2 
    \end{array} \hspace{4cm}
        \begin{array}{c|cc}
        \gamma & \gamma & 0  \\
        1 & 1-\gamma & \gamma \\ \hline
         & 1-\gamma & \gamma 
    \end{array}
\end{equation}
and taking $\gamma = (2-\sqrt{2})/2$ guarantees that the implicit part of the scheme is L-stable.
\subsubsection*{I3: CK(3,3,2)}
Second order L-stable GSA IMEX-RK type II, with $\gamma = (2-\sqrt{2})/2$.
\begin{equation}
    \begin{array}{c|ccc}
        0 & 0 & 0 & 0 \\
        2/3 & 2/3 & 0 & 0 \\
        1 & 1/4 & 3/4 & 0  \\ \hline
         & 1/4 & 3/4 & 0
    \end{array} \hspace{4cm}
        \begin{array}{c|ccc}
        0 & 0 & 0 & 0 \\
        2/3 & 2/3-\gamma & \gamma & 0  \\
        1 & 1/4-\gamma/2 & 3/4-3/2\gamma & \gamma \\ \hline
         & 1/4-\gamma/2 & 3/4-3/2\gamma & \gamma 
    \end{array}
\end{equation}

\subsubsection*{I4: IMEX-SSP2(3,3,2)}
Second order GSA IMEX-RK method of type I.
\begin{equation}
    \begin{array}{c|ccc}
        0 & 0 & 0 & 0 \\
        1/2 & 1/2 & 0 & 0 \\
        1 & 1/2 & 1/2 & 0  \\ \hline
         & 1/3 & 1/3 & 1/3
    \end{array} \hspace{4cm}
        \begin{array}{c|ccc}
        1/4 & 1/4 & 0 & 0 \\
        1/4 & 0 & 1/4 & 0  \\
        1 & 1/3 & 1/3 & 1/3  \\ \hline
         & 1/3 & 1/3 & 1/3 
    \end{array}
\end{equation}

\subsubsection*{I5: IMEX-SSP3(4,3,3)}
Third order IMEX-RK method of type I. This method is called IMEX-SSP3(4,3,3) where the explicit method is an SSP-RK
\begin{equation}
        \begin{array}{c|cccc}
        0 & 0 & 0 & 0  & 0\\
        0 & 0 & 0 & 0  & 0 \\
         1 & 0 & 1 & 0 & 0 \\ 
         1/2 & 0 & 1/4 & 1/4 & 0\\ \hline
          & 0 & 1/6 & 1/6 & 2/3
    \end{array} \hspace{2cm}
        \begin{array}{c|cccc}
        \alpha & \alpha & 0 & 0  & 0\\
        0 & -\alpha & \alpha & 0 & 0 \\
        1 & 0 & 1-\alpha & \alpha & 0 \\
        1/2 & \beta & \eta & 1/2-\alpha-\beta-\eta & \alpha \\ \hline
         & 0 & 1/6 & 1/6 & 2/3
    \end{array}
\end{equation}
with $\alpha = 0.24169426078821, \beta = 0.06042356519705, \eta = 0.12915286960590.$

\subsubsection*{I6: SI-IMEX(4,4,3)}
Third order L-stable SA IMEX-RK scheme of type I
\begin{equation}
\begin{array}{c|cccc}
        0 & 0 & 0 & 0  & 0\\
        \gamma & \gamma & 0 & 0  & 0 \\
        \beta & \delta & \xi & 0 & 0 \\
        1 & \zeta & \tau & \theta & 0 \\ \hline
         & 0 & \alpha & \eta & \gamma
    \end{array} \hspace{4cm}
        \begin{array}{c|cccc}
        \gamma & \gamma & 0 & 0  & 0\\
        \gamma & 0 & \gamma & 0  & 0 \\
        \beta & 0 & \mu & \gamma & 0\\ 
         1 & 0 & \alpha & \eta & \gamma \\ \hline
         & 0 & \alpha & \eta & \gamma
    \end{array}
\end{equation}
$\alpha = 1.208496649176,\beta = 0.717933260754, \delta = 1.243893189, \eta = -0.644363170684, \gamma = 0.435866521508, \mu = 0.282066739245, \xi = -0.5259599287, \zeta = 0.6304125582, \tau = 0.7865807402, \theta = -0.4169932983.$

\subsubsection{Standard time discretization}
We now briefly present a standard Implicit Euler discretization of system~\eqref{eq:FEM_semidiscr_1}. In this case, we compute the three unknowns of the system in separate stages.

Given $c_{\pm,h}^{n}$ we find $c_{\pm,h}^{n+1}$ as follows: 
\begin{itemize} 
\item compute the Coulomb potential $\Phi^n_h$ by solving the discrete Poisson equation 
\begin{equation}
\label{eq:serial}
	 \varepsilon\mathbb L[u_h]\, \Phi^n_h = \mathbb B[u_h] \, \left( \frac{c_{+,h}}{m_+} - \frac{c_{-,h}}{m_-}\right)
\end{equation}
\item compute the concentrations $c_{\pm,h}^{n+1}$ by solving the drift-diffusion equations
\begin{align}
\label{eq:serial1}
\mathbb B[v_h] c^{n+1}_{\pm,h} = &\mathbb B[v_h] c^{n}_{\pm,h} - \frac 1 2\Delta t D_\pm \left( \mathbb L[v_h]  \pm \mathbb G\left[ \Phi^n_h ,v_h \right]\right)  c^{n+1}_{\pm,h} \\\nonumber & -\frac 1 2\Delta t D_\pm \left( \mathbb L[v_h]  \pm \mathbb G\left[ \Phi^n_h ,v_h \right]\right)  c^{n}_{\pm,h} 
\end{align}
\end{itemize}
where the operator $\mathbb G[\Phi_h,v_h]$ is the operator such that \\$\mathbb G \left[ \Phi_h,v_h \right] c_{\pm,h} \approx \left(c_{\pm,h} \nabla \Phi_h,\nabla v_h  \right)_{L^2\left(\Omega_h\right)}$. This time we consider the concentrations as implicit part of the drift term.
\begin{figure}[h]
\centering
\includegraphics[width=0.55\textwidth]{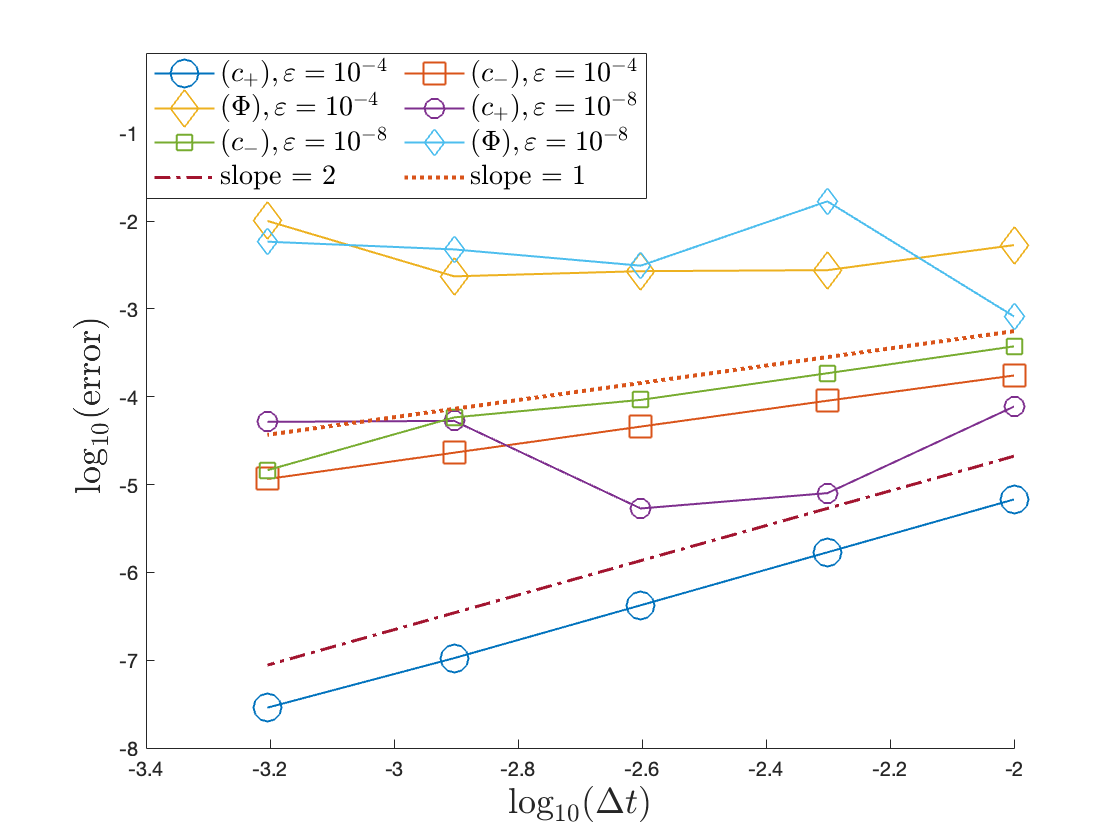}
\caption{\textit{Time accuracy orders of the system~(\ref{eq:serial}-\ref{eq:serial1}) at final time $t = 0.1$, for different values of $\varepsilon$.}}  
\label{fig:accuracycpm_seriale}
\end{figure}
\section{Numerical Results}
\label{sec:numericalres}
In this section, we compare the accuracy order of the three time discretizations that we describe in the previous section. Our goal is to show which numerical schemes are stable for vanishing values of $\varepsilon$ (AP property) and which ones maintain the same accuracy order in the same limit (AA property). 

The initial conditions for Eqs.~(\ref{eq:vectorial}-\ref{eq:vectorial1}) and for Eqs.~(\ref{eq:serial}-\ref{eq:serial1}) are given by
\begin{align} 	
\label{eq_initial_2D}
c_\pm^{\rm in}(x,y,t=0) &= \frac{v_0}{2\sigma^2}\exp\left(-{\left((x-x^{\rm in}_\pm)^2 + (y-y^{\rm in}_\pm)^2 \right)}/{2\sigma^2}\right), 
\end{align}
while, for Eqs.~\eqref{eq:vectorial},\eqref{eq:vectorial2} are
\begin{subequations}
    \begin{align}
        \mathcal{C}^{\rm in}(x,y,t=0) &= \frac{c_+^{\rm in}(x,y,t=0)}{m_+} + \frac{c_-^{\rm in}(x,y,t=0)}{m_-} \\
        \mathcal{Q}^{\rm in}(x,y,t=0) &= \frac{1}{\varepsilon} \left( \frac{c_+^{\rm in}(x,y,t=0)}{m_+} - \frac{c_-^{\rm in}(x,y,t=0)}{m_-} \right)
    \end{align}
\end{subequations}
where $v_0$ denotes the total volume per unit surface. In our numerical tests $v_0 = 10^{-6}, \, x_\pm^{\rm in} = 0.5 \mp 0.1,\, y_\pm^{\rm in}  = 0.2,\, \sigma = 0.05$, the number of cells of the space discretization is $N = 100$ and $\Delta t = h$, unless otherwise specified. The other parameters are expressed in Table~\ref{table_parameters}.

\begin{figure}[h]
\centering
\begin{minipage}{.49\textwidth}
\begin{overpic}[abs,width=\textwidth,unit=1mm,scale=.25]{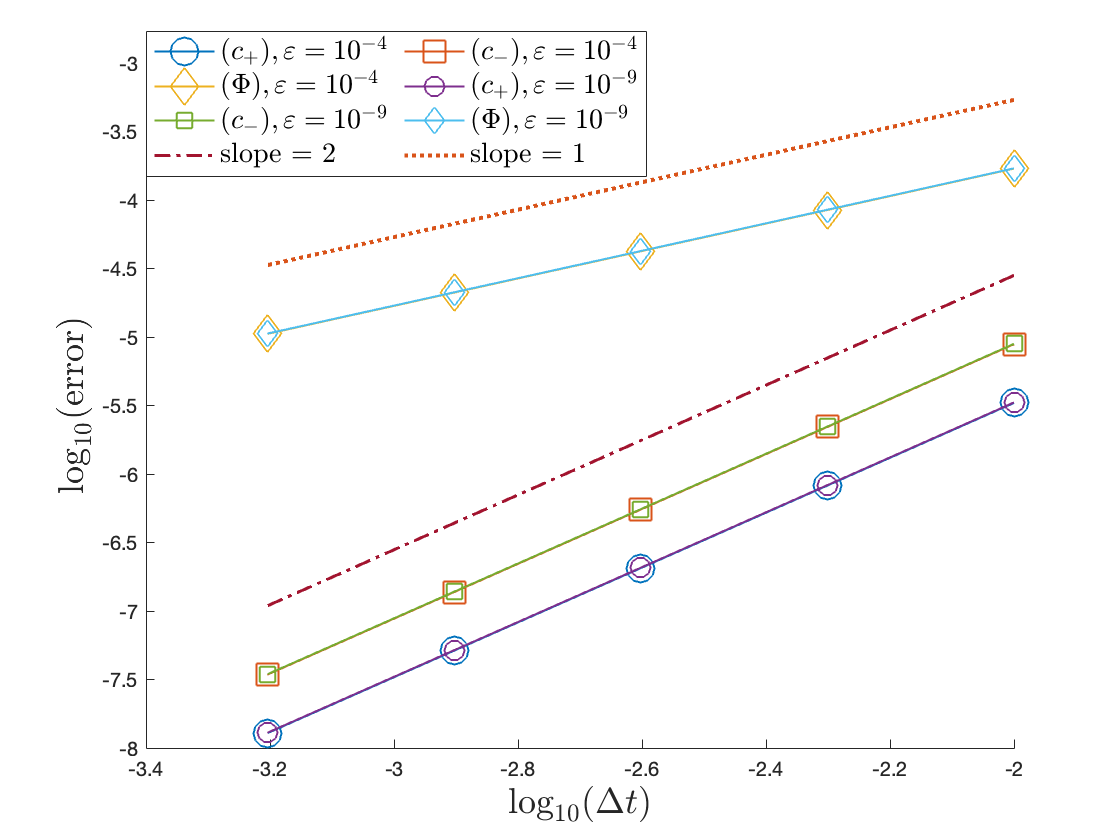}  
\put(-1,42){(a)}
\end{overpic}
    \end{minipage}
\begin{minipage}{.49\textwidth}
\begin{overpic}[abs,width=\textwidth,unit=1mm,scale=.25]{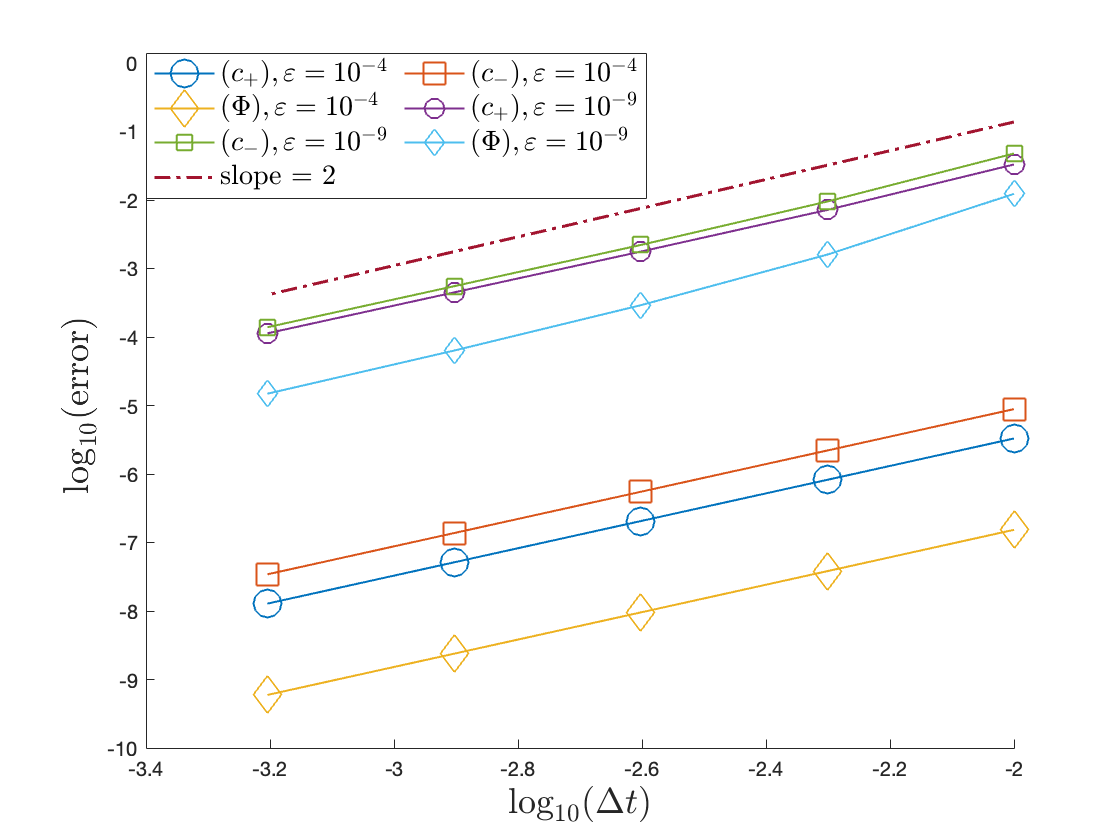}  
\put(-1,42){(b)}
\end{overpic}
    \end{minipage}
    \begin{minipage}{.49\textwidth}
\begin{overpic}[abs,width=\textwidth,unit=1mm,scale=.25]{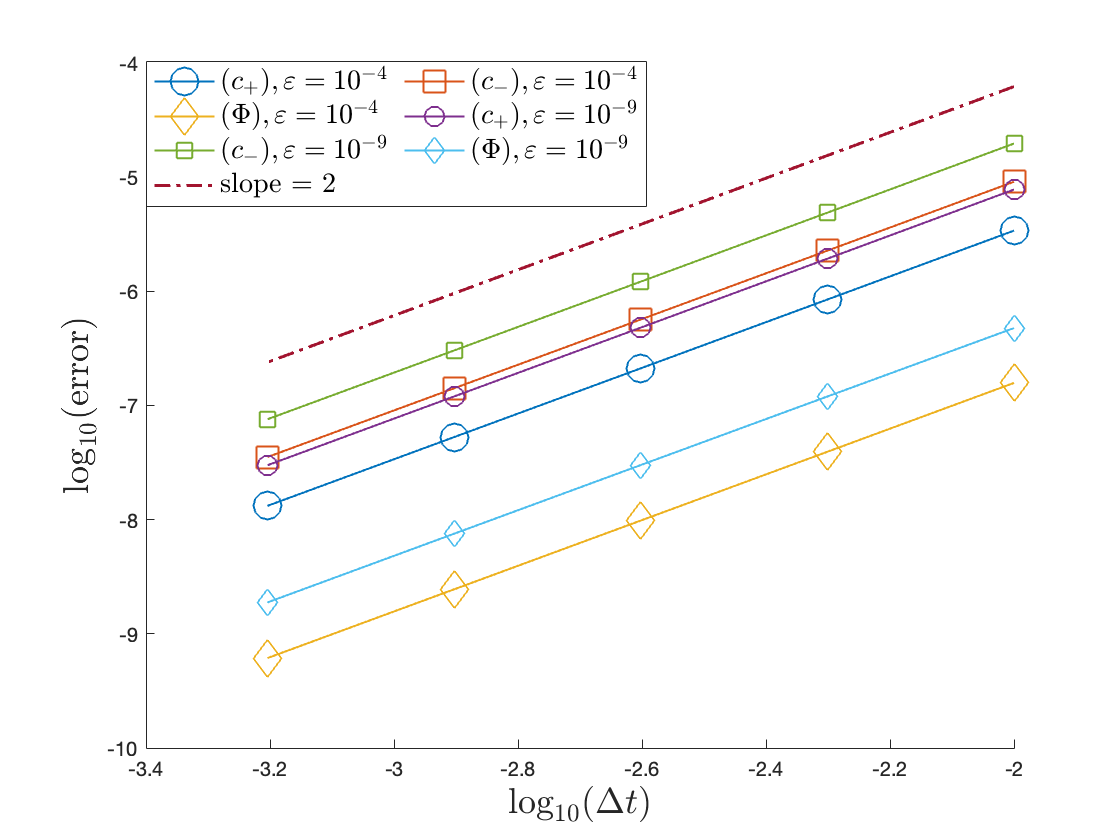}     
\put(-1,42){(c)}
\end{overpic}
    \end{minipage}
    \begin{minipage}{.49\textwidth}
\begin{overpic}[abs,width=\textwidth,unit=1mm,scale=.25]{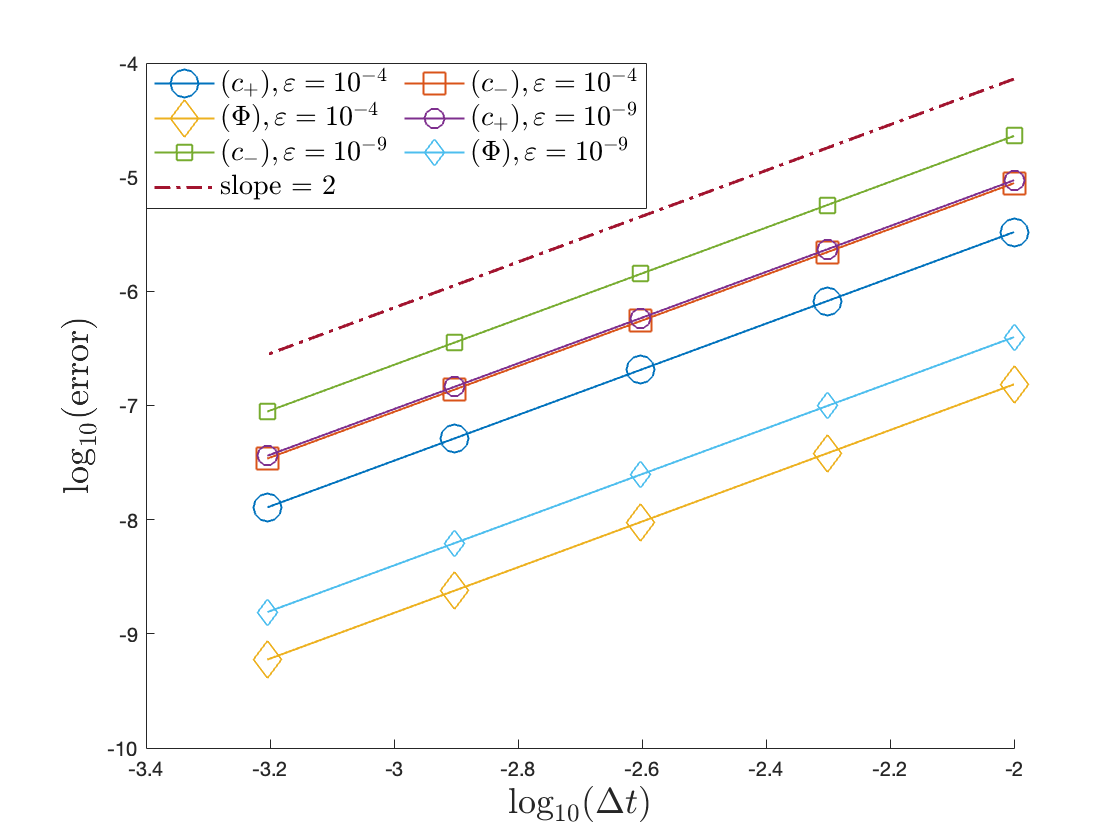}     
\put(-1,42){(d)}
\end{overpic}
    \end{minipage}
\begin{minipage}{.49\textwidth}
\begin{overpic}[abs,width=\textwidth,unit=1mm,scale=.25]{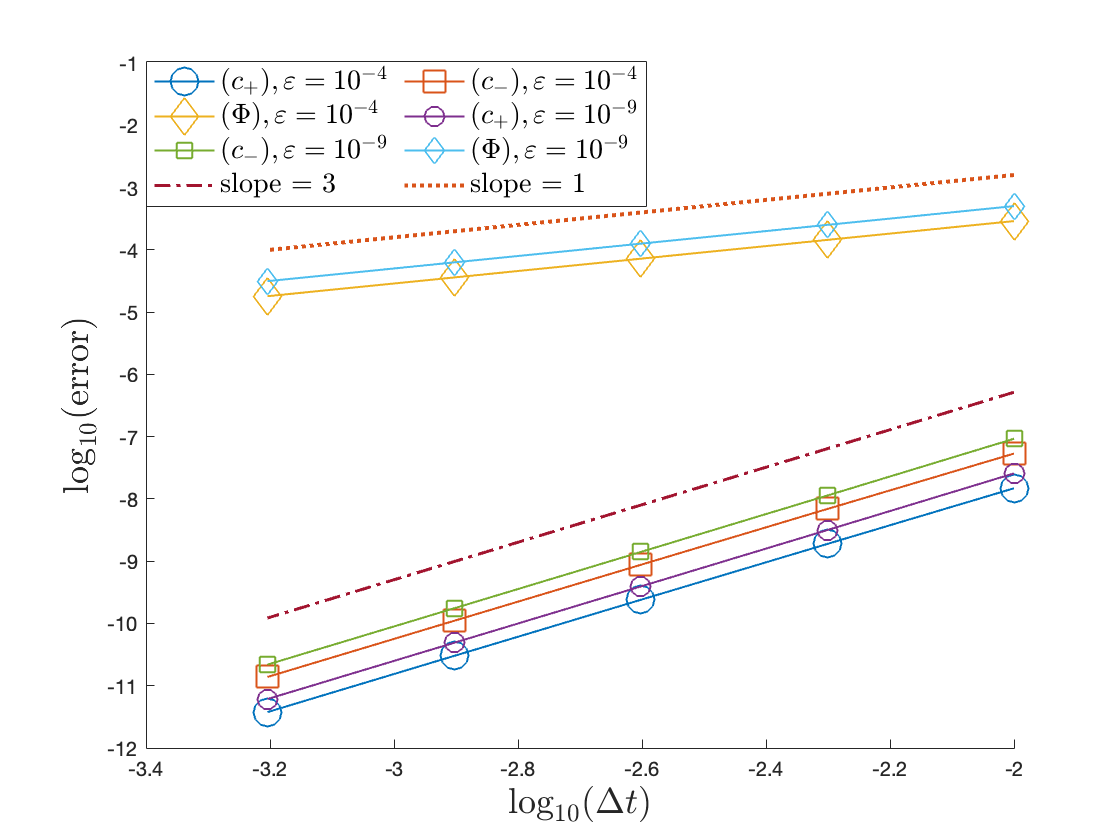}  
\put(-1,42){(e)}
\end{overpic}
    \end{minipage}
    \begin{minipage}{.49\textwidth}
\begin{overpic}[abs,width=\textwidth,unit=1mm,scale=.25]{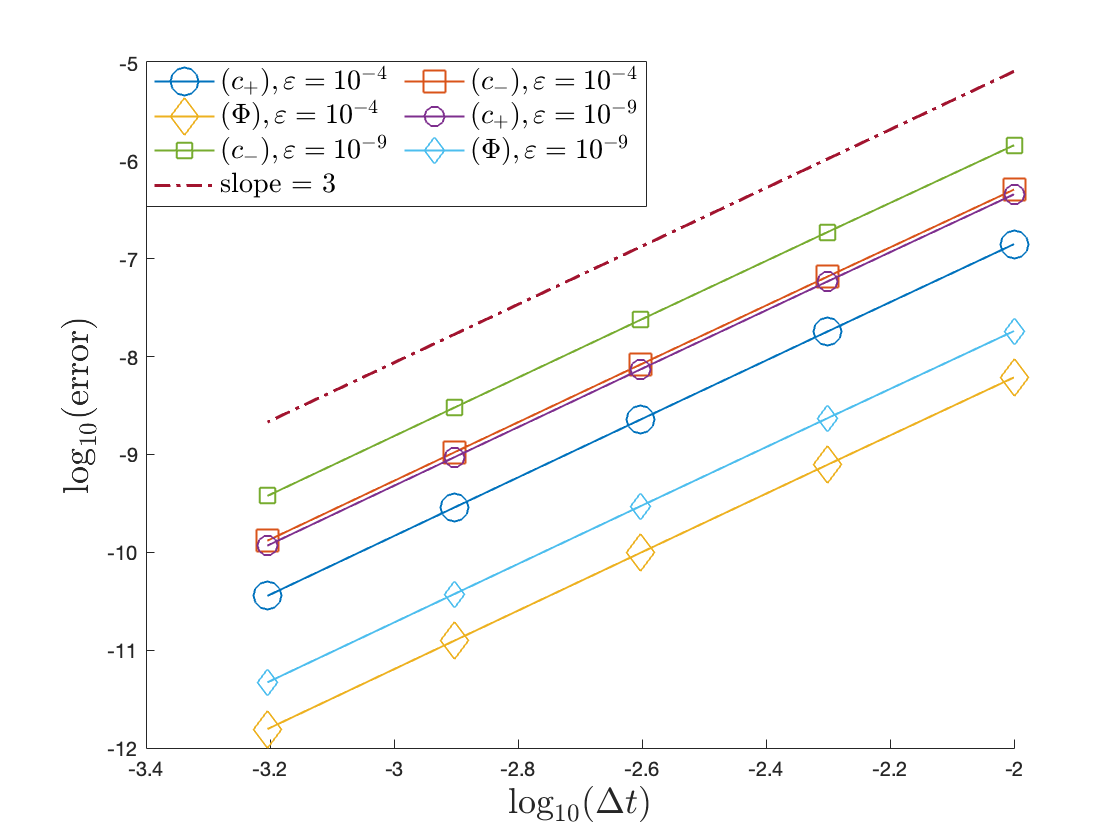}     
\put(-1,42){(f)}
\end{overpic}
    \end{minipage}
\caption{\textit{Time accuracy orders of the system~(\ref{eq:vectorial}-\ref{eq:vectorial1}) at final time $t = 0.1$, for different values of $\varepsilon$. }}  
\label{fig:accuracy_cpm}
\end{figure}

\begin{figure}[h]
\centering
\begin{minipage}{.49\textwidth}
\begin{overpic}[abs,width=\textwidth,unit=1mm,scale=.25]{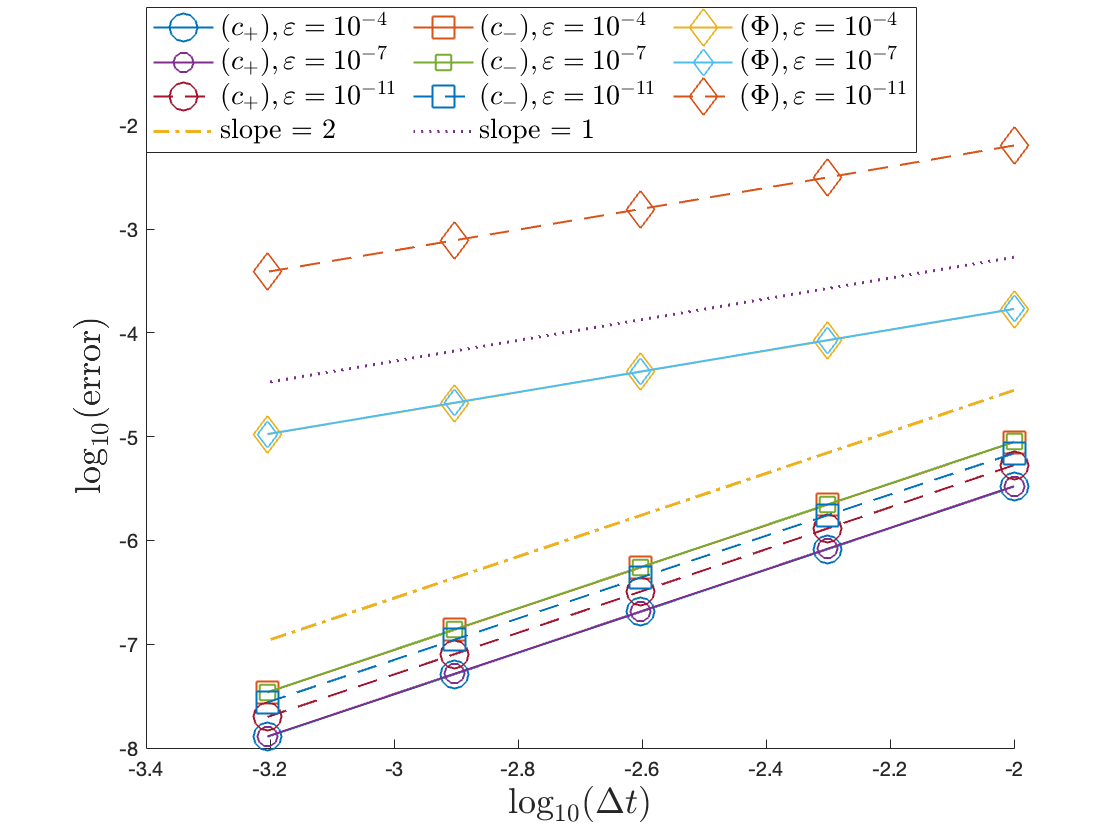}  
\put(-1,42){(a)}
\end{overpic}
    \end{minipage}
\begin{minipage}{.49\textwidth}
\begin{overpic}[abs,width=\textwidth,unit=1mm,scale=.25]{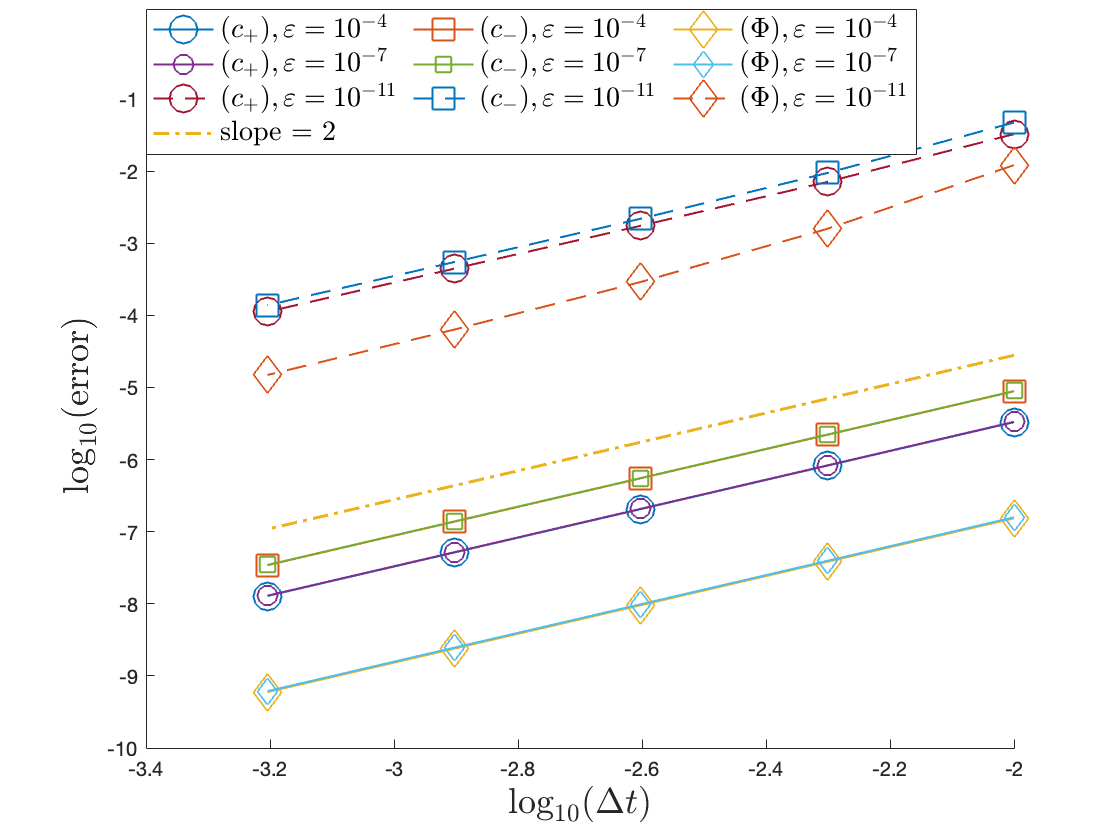}  
\put(-1,42){(b)}
\end{overpic}
    \end{minipage}
    \begin{minipage}{.49\textwidth}
\begin{overpic}[abs,width=\textwidth,unit=1mm,scale=.25]{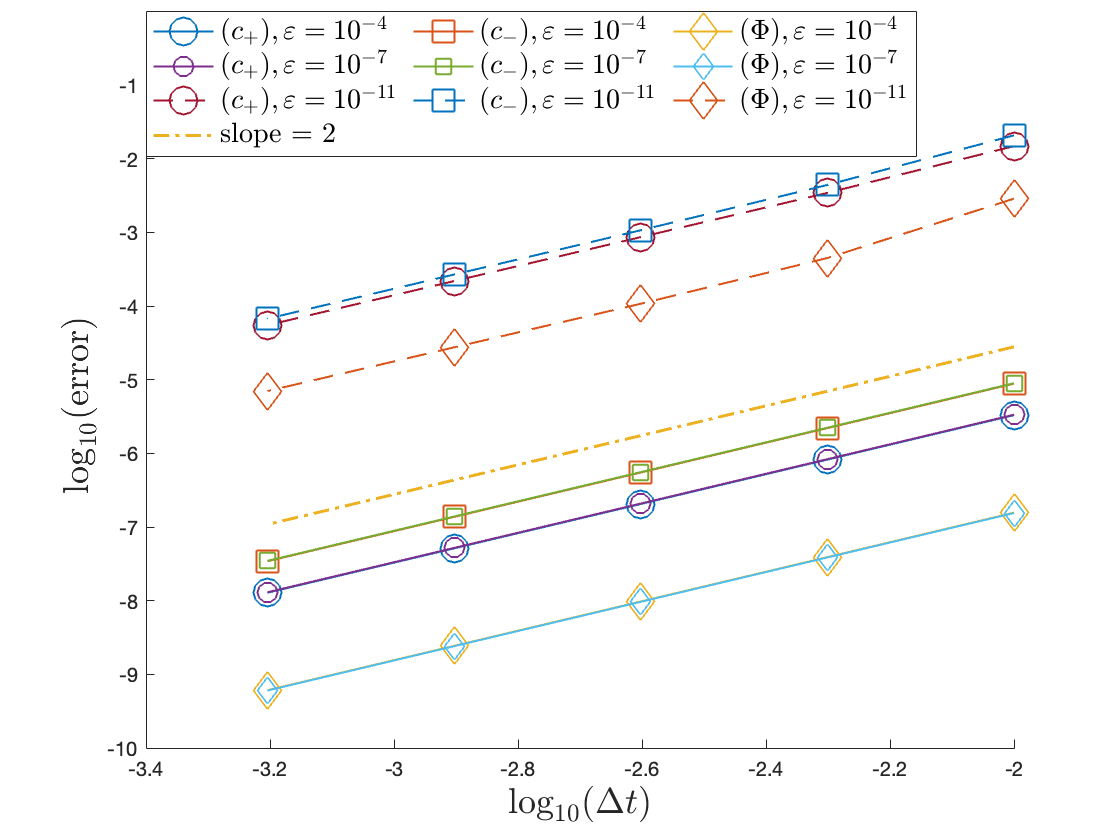}     
\put(-1,42){(c)}
\end{overpic}
    \end{minipage}
    \begin{minipage}{.49\textwidth}
\begin{overpic}[abs,width=\textwidth,unit=1mm,scale=.25]{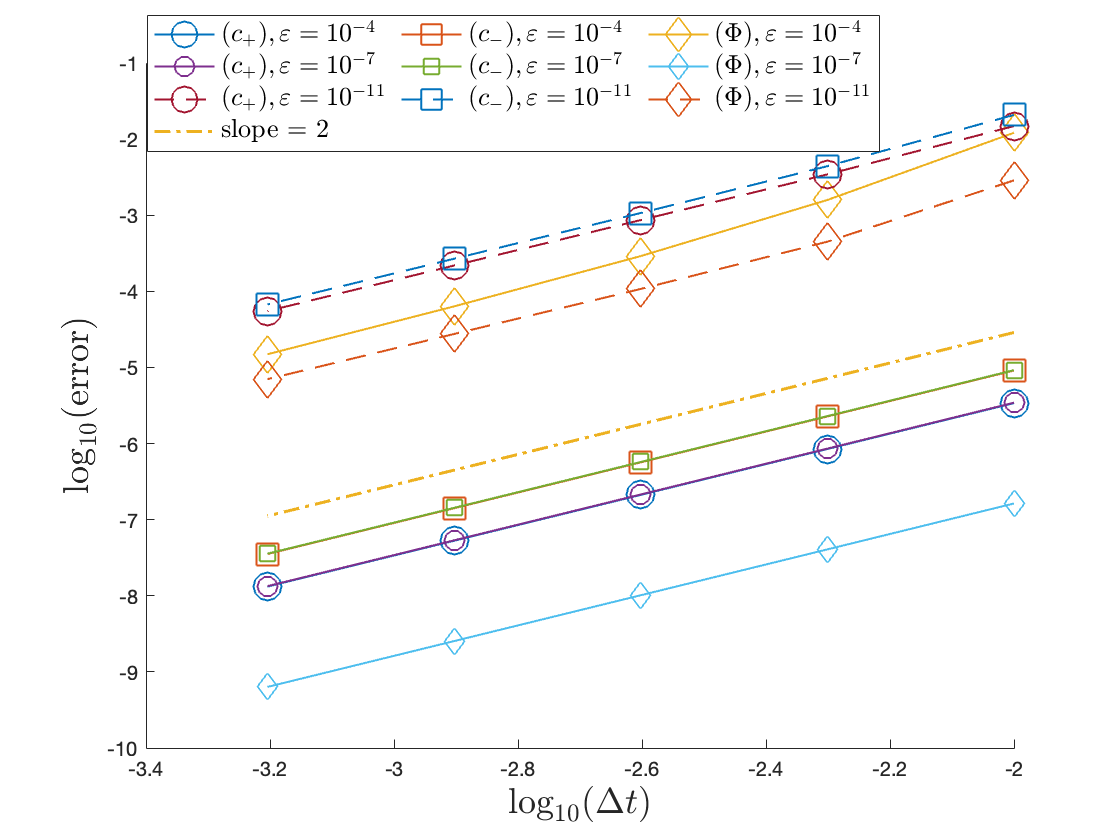}     
\put(-1,42){(d)}
\end{overpic}
    \end{minipage}
\begin{minipage}{.49\textwidth}
\begin{overpic}[abs,width=\textwidth,unit=1mm,scale=.25]{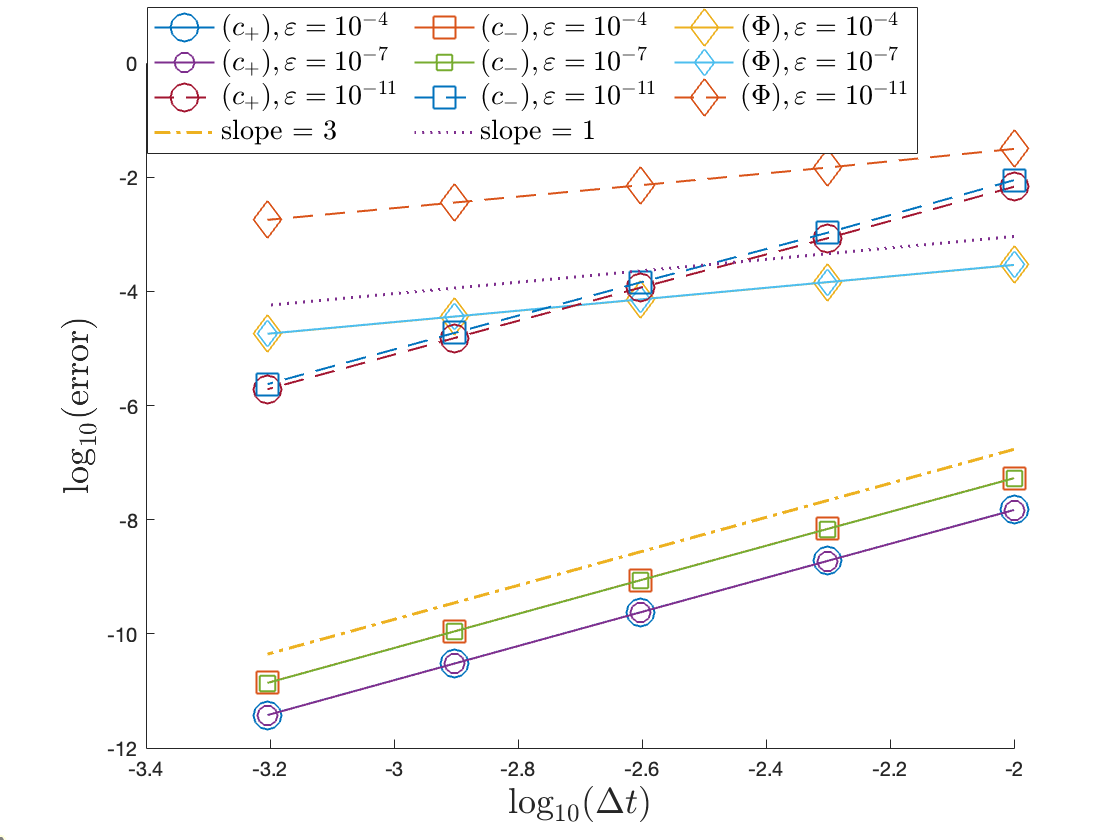}  
\put(-1,42){(e)}
\end{overpic}
    \end{minipage}
    \begin{minipage}{.49\textwidth}
\begin{overpic}[abs,width=\textwidth,unit=1mm,scale=.25]{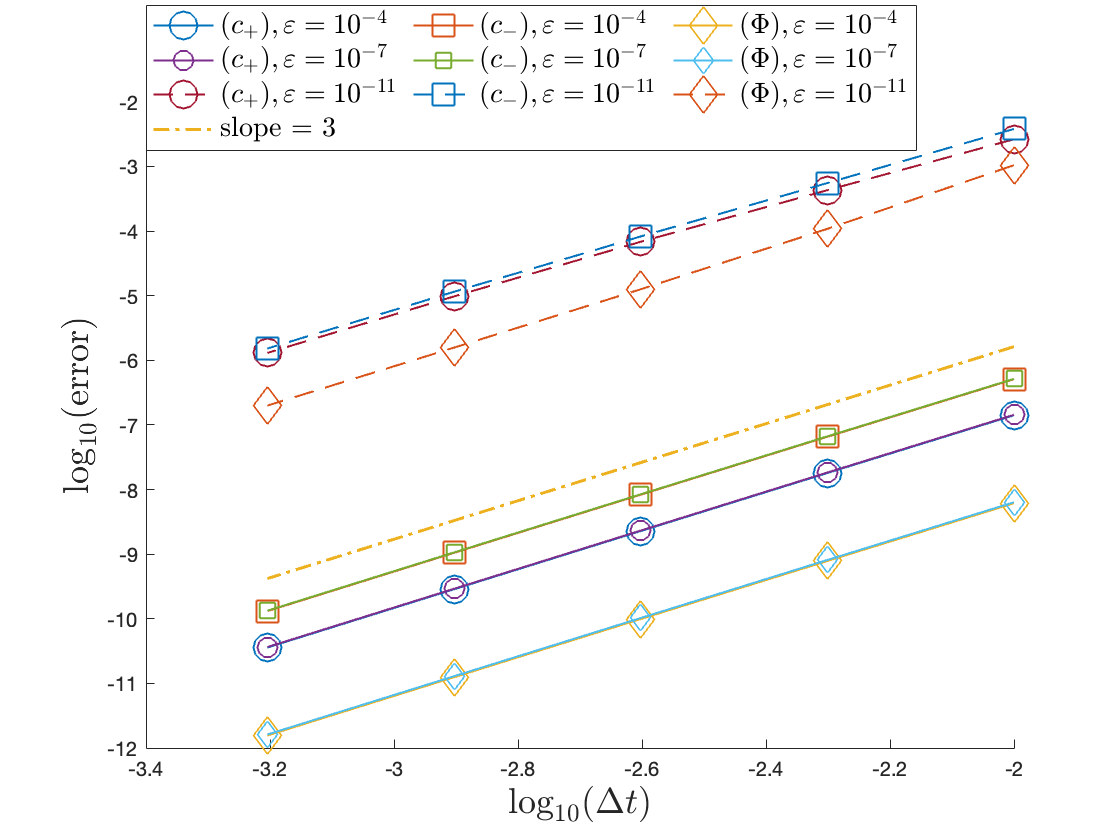}     
\put(-1,42){(f)}
\end{overpic}
    \end{minipage}
\caption{\textit{Time accuracy orders of the system~\eqref{eq:vectorial},\eqref{eq:vectorial2} at final time $t = 0.1$, for different values of $\varepsilon$. }}  
\label{fig:accuracy}
\end{figure}

In Figs.~(\ref{fig:accuracycpm_seriale}-\ref{fig:accuracy}) we show the $L^2$-norm of the error as a function of $\Delta t$ (and fixed $h = 10^{-2}$), at final time $t = 0.1$ and for different values of $\varepsilon$. In the absence of an exact solution, we apply the Richardson extrapolation technique (see, e.g., \cite{richardson1911ix}) to estimate the order of the method. We observe that the numerical method in Eqs.~\eqref{eq:vectorial},\eqref{eq:vectorial2} remains stable (AP) and is second-order accurate (AA) for all the tested values of $\varepsilon$. In contrast, the scheme in Eqs.~(\ref{eq:vectorial}-\ref{eq:vectorial1})  loses stability for $\varepsilon < 10^{-9}$, although it still achieves second-order accuracy within the considered range of $\varepsilon$. The scheme in Eqs.~(\ref{eq:serial}-\ref{eq:serial1}) becomes unstable for $\varepsilon < 10^{-8}$ and additionally exhibits a deterioration in its order of accuracy for the tested values of $\varepsilon$.
Thus, we conclude that the formulation in Eqs.~\eqref{eq:vectorial},\eqref{eq:vectorial2} is the best choice for stability and accuracy.

Table~\ref{tab:time} reports the execution time per iteration for the two systems (\ref{eq:vectorial}-\ref{eq:vectorial1}) (denoted as $c_\pm$-system) and \eqref{eq:vectorial},\eqref{eq:vectorial2} (denoted as $\mathcal C, \mathcal Q$-system) as the parameter $\varepsilon$ varies. It highlights how the computational cost increases as $\varepsilon$ decreases, showing that the original formulation of the Poisson-Nernst-Planck system is more efficient for non-negligible values of the Debye length. We limit our analysis to $\varepsilon \ge 10^{-9}$, since for smaller values the time step restriction for the $c_\pm$-model becomes prohibitively small ($\Delta t = 10^{-4}h$).

\begin{figure}[h]
\centering
\begin{minipage}{.49\textwidth}
\begin{overpic}[abs,width=\textwidth,unit=1mm,scale=.25]{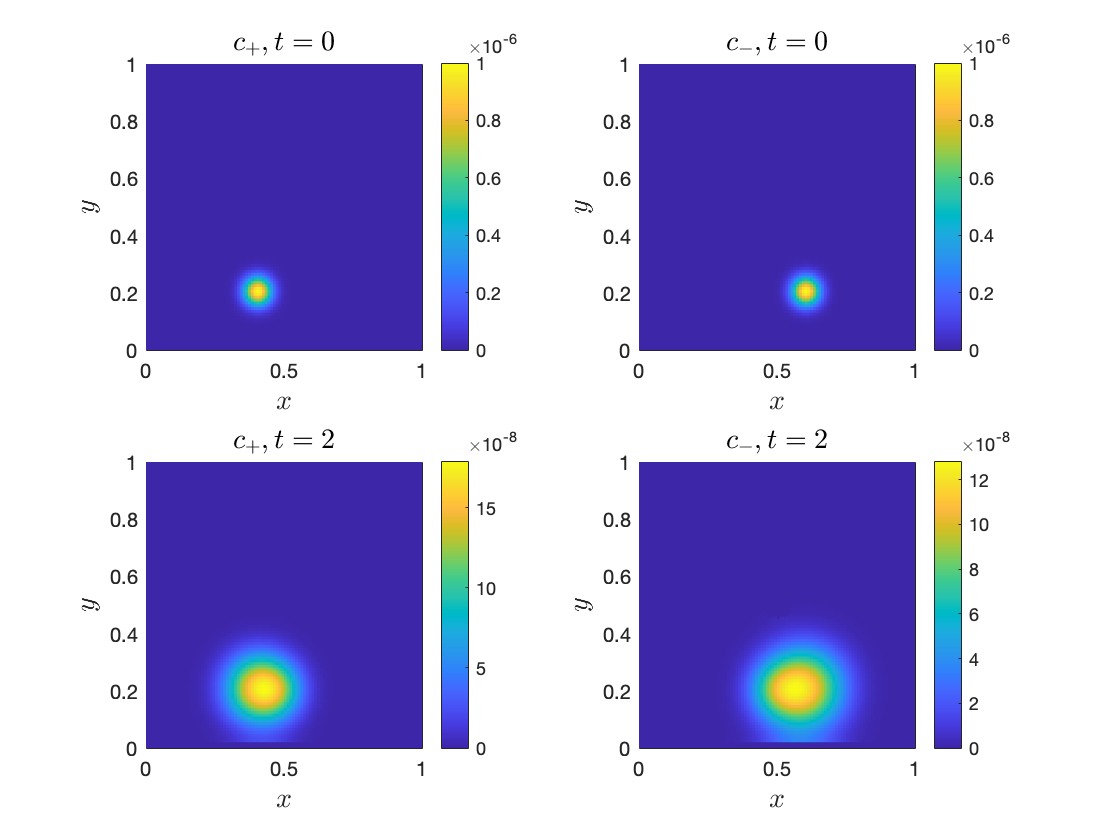}  
\end{overpic}
    \end{minipage}
\begin{minipage}{.49\textwidth}
\begin{overpic}[abs,width=\textwidth,unit=1mm,scale=.25]{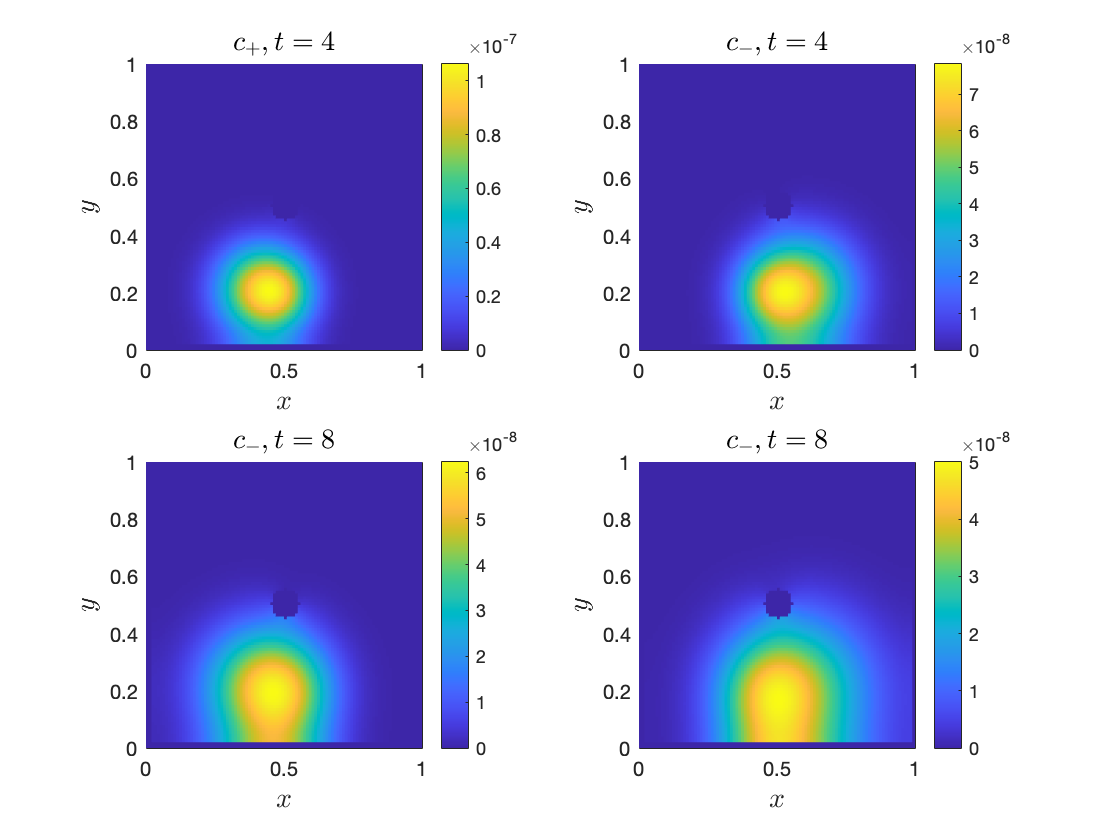}  
\end{overpic}
    \end{minipage}
\caption{\textit{Time evolution of the system~\eqref{eq:vectorial},\eqref{eq:vectorial2} reaching the quasi-neutrality limit, with $\varepsilon = 10^{-9}$. Time discretization chosen is I2: IMEX-SA(2,2,2).}}  
\label{fig:solutimes9}
\end{figure}

\begin{figure}[h]
\centering
\begin{minipage}{.49\textwidth}
\begin{overpic}[abs,width=\textwidth,unit=1mm,scale=.25]{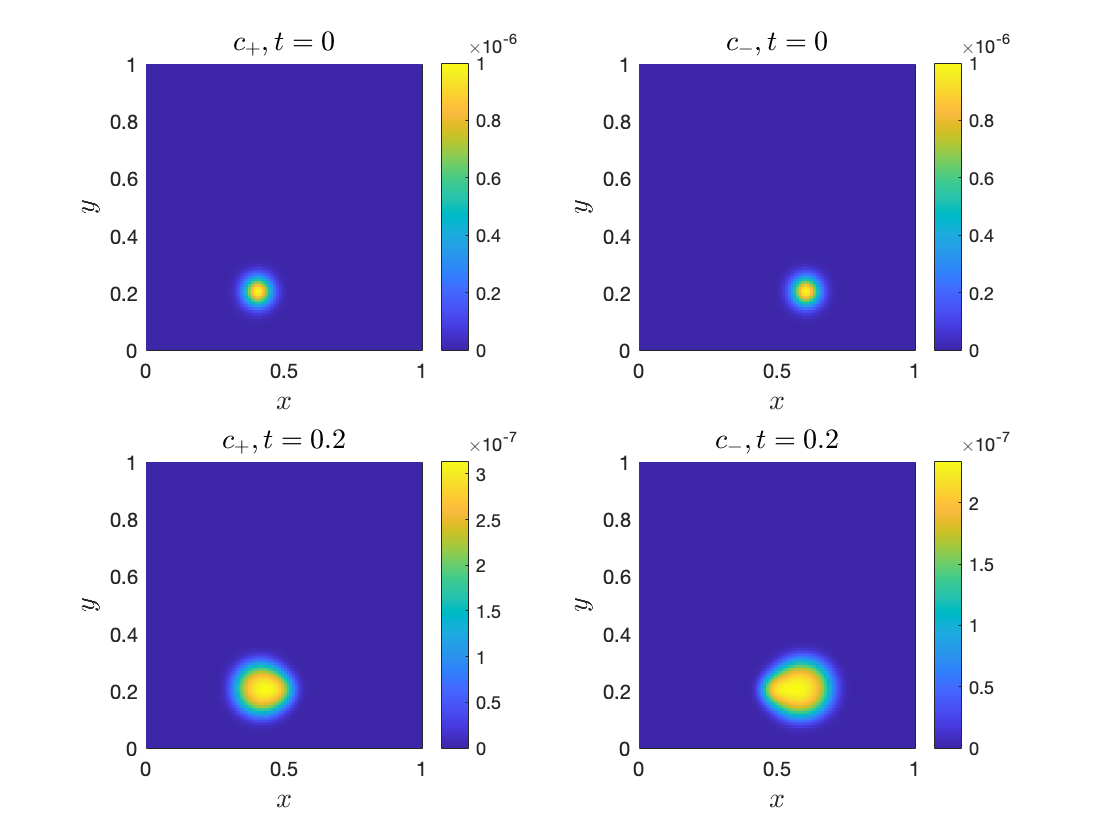}  
\end{overpic}
    \end{minipage}
\begin{minipage}{.49\textwidth}
\begin{overpic}[abs,width=\textwidth,unit=1mm,scale=.25]{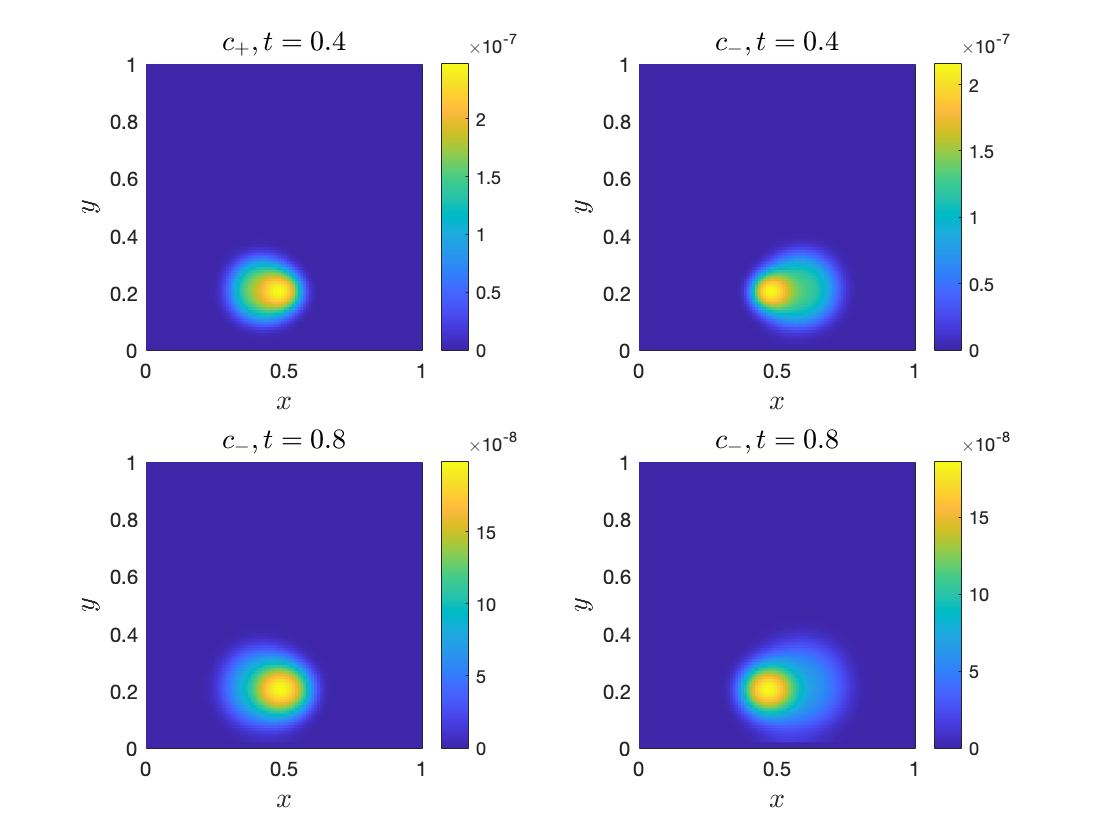}  
\end{overpic}
    \end{minipage}
\caption{\textit{Time evolution of the system~\eqref{eq:vectorial},\eqref{eq:vectorial2} reaching the quasi-neutrality limit, with $\varepsilon = 10^{-10}$. Time discretization chosen is I2: IMEX-SA(2,2,2).}}  
\label{fig:solutimes10}
\end{figure}

\begin{figure}[h]
\centering
\begin{minipage}{.49\textwidth}
\begin{overpic}[abs,width=\textwidth,unit=1mm,scale=.25]{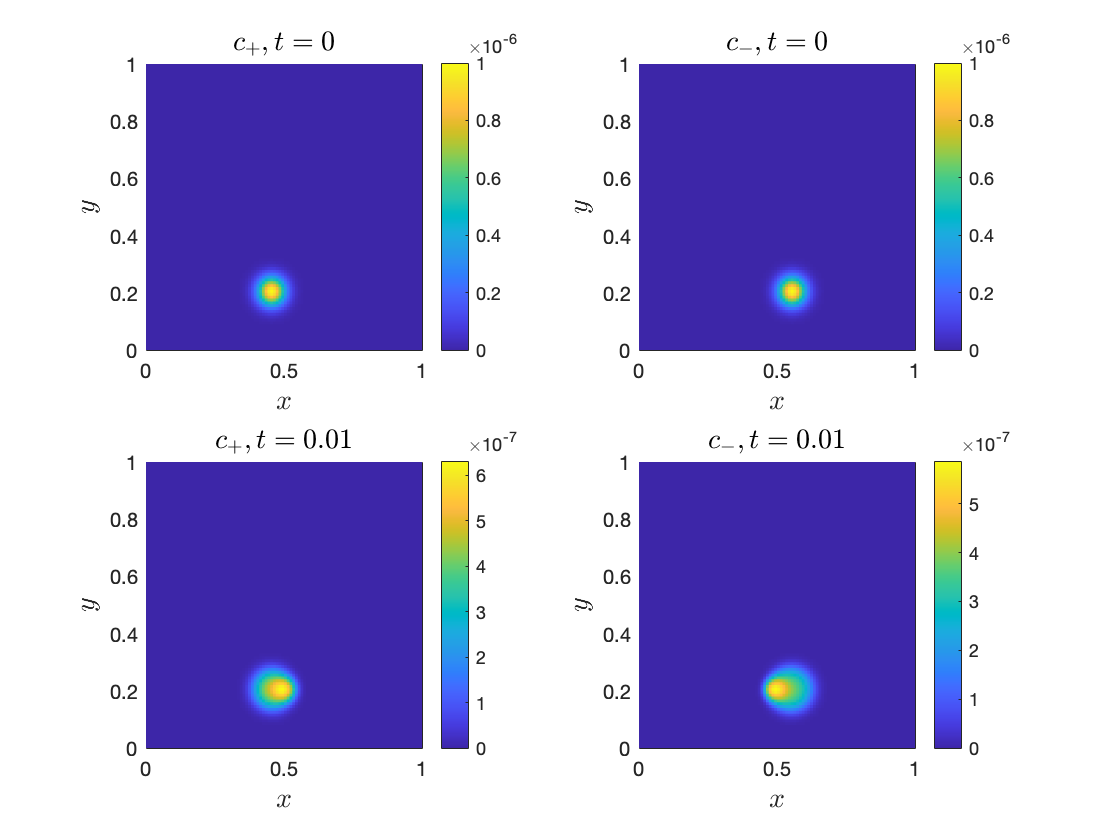}  
\end{overpic}
    \end{minipage}
\begin{minipage}{.49\textwidth}
\begin{overpic}[abs,width=\textwidth,unit=1mm,scale=.25]{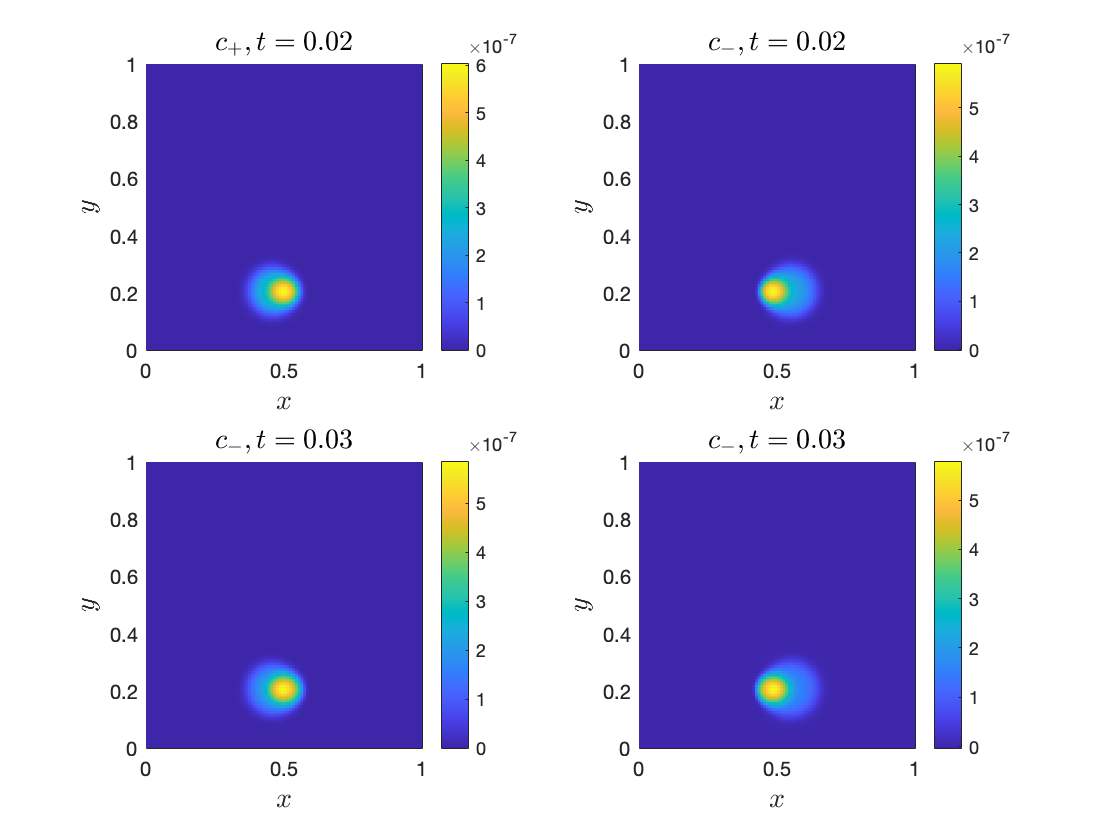}  
\end{overpic}
    \end{minipage}
\caption{\textit{Time evolution of the system~\eqref{eq:vectorial},\eqref{eq:vectorial2} reaching the quasi-neutrality limit, with $\varepsilon = 10^{-11}$. Time discretization chosen is I2: IMEX-SA(2,2,2), $x^{\rm in}_\pm = 0.5 \mp 0.05$ and $\Delta t = 0.1h$.}}  
\label{fig:solutimes11_3}
\end{figure}

\begin{figure}[h]
\centering
\begin{minipage}{.49\textwidth}
\begin{overpic}[abs,width=\textwidth,unit=1mm,scale=.25]{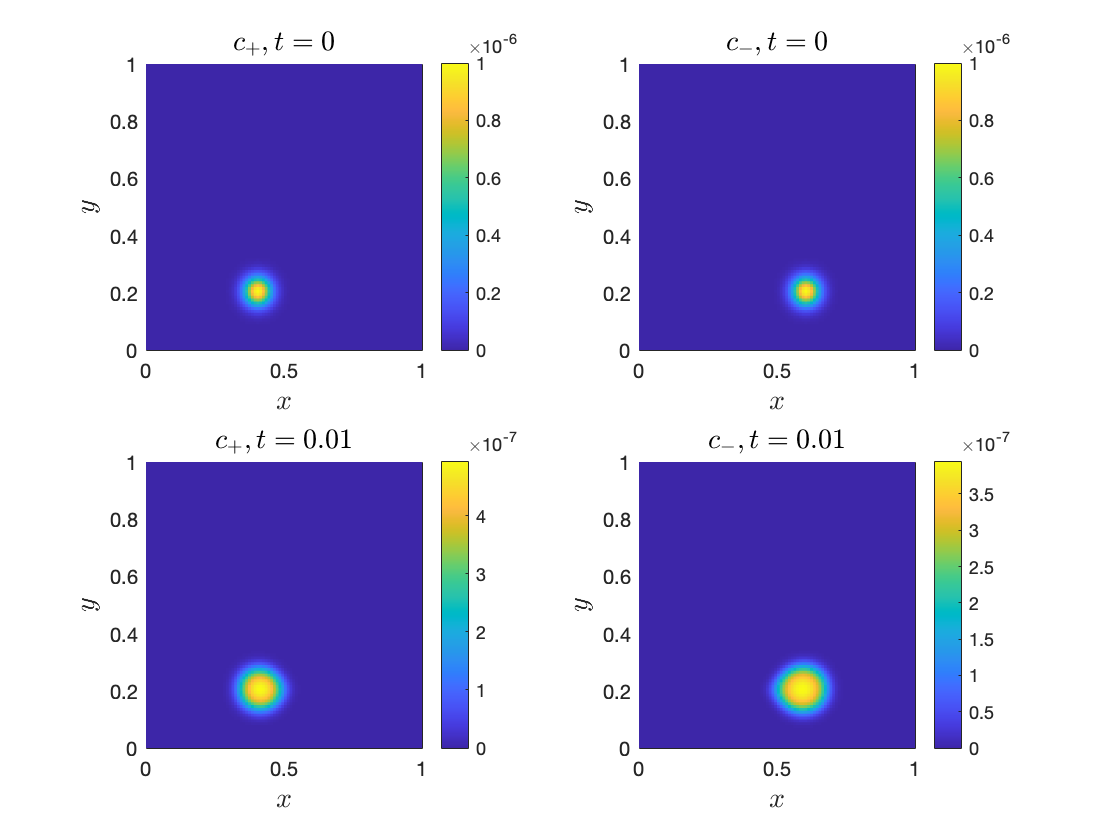}  
\end{overpic}
    \end{minipage}
\begin{minipage}{.49\textwidth}
\begin{overpic}[abs,width=\textwidth,unit=1mm,scale=.25]{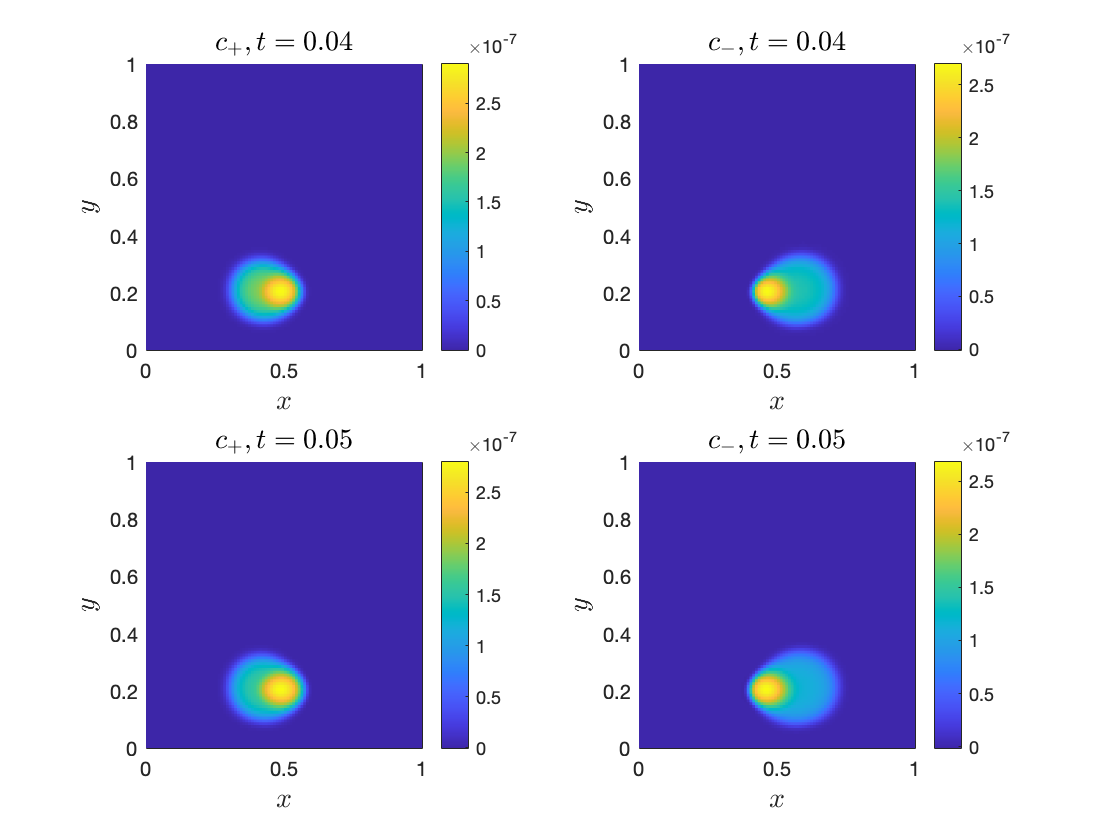}  
\end{overpic}
    \end{minipage}
\begin{minipage}{.49\textwidth}
\begin{overpic}[abs,width=\textwidth,unit=1mm,scale=.25]{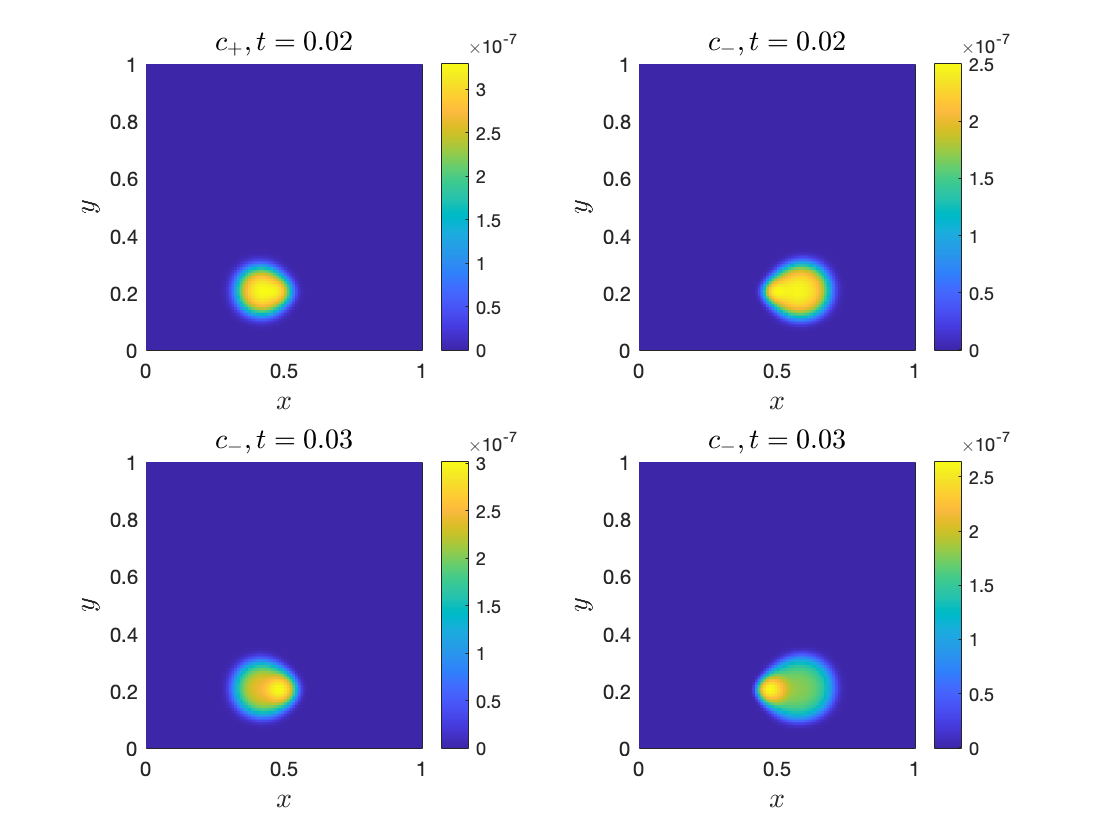}  
\end{overpic}
    \end{minipage}
\begin{minipage}{.49\textwidth}
\begin{overpic}[abs,width=\textwidth,unit=1mm,scale=.25]{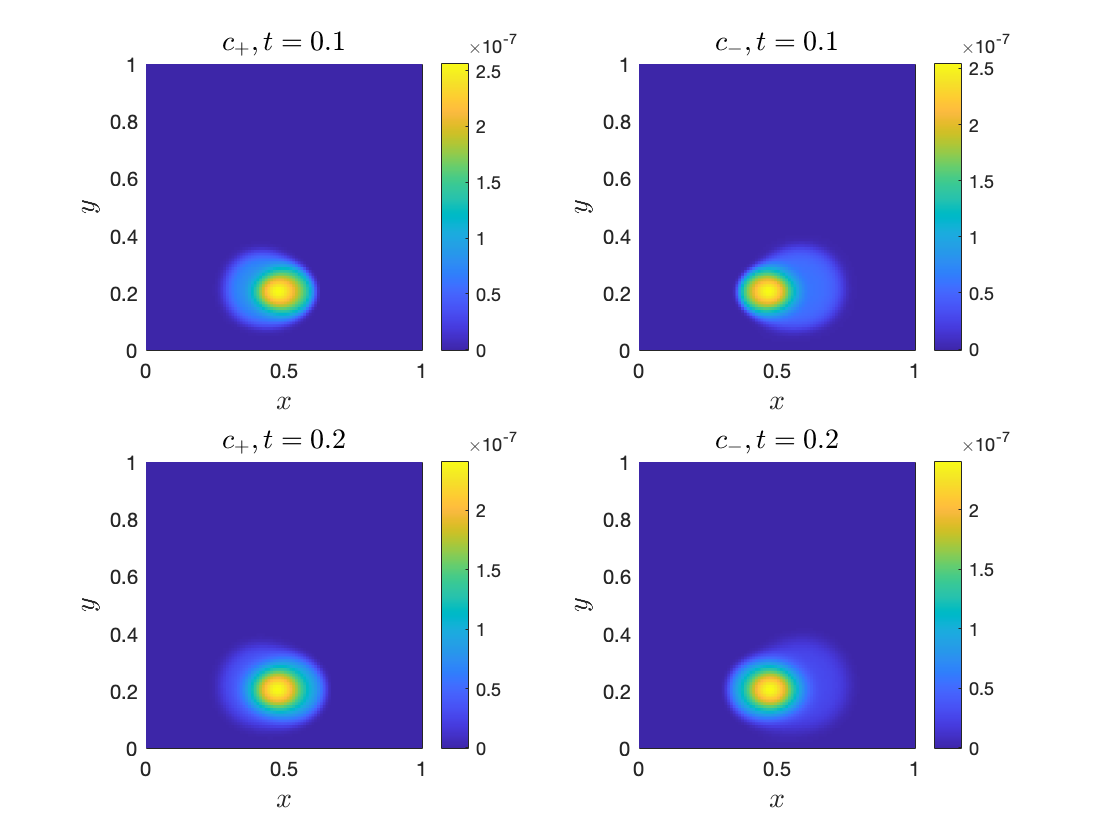}  
\end{overpic}
    \end{minipage}    
\caption{\textit{Time evolution of the system~\eqref{eq:vectorial},\eqref{eq:vectorial2} reaching the quasi-neutrality limit, with $\varepsilon = 10^{-11}$. Time discretization chosen is I2: IMEX-SA(2,2,2) and $\Delta t = 0.1h$.}}  
\label{fig:solutimes11}
\end{figure}

\begin{figure}[h]
\centering
\begin{minipage}{.49\textwidth}
\begin{overpic}[abs,width=\textwidth,unit=1mm,scale=.25]{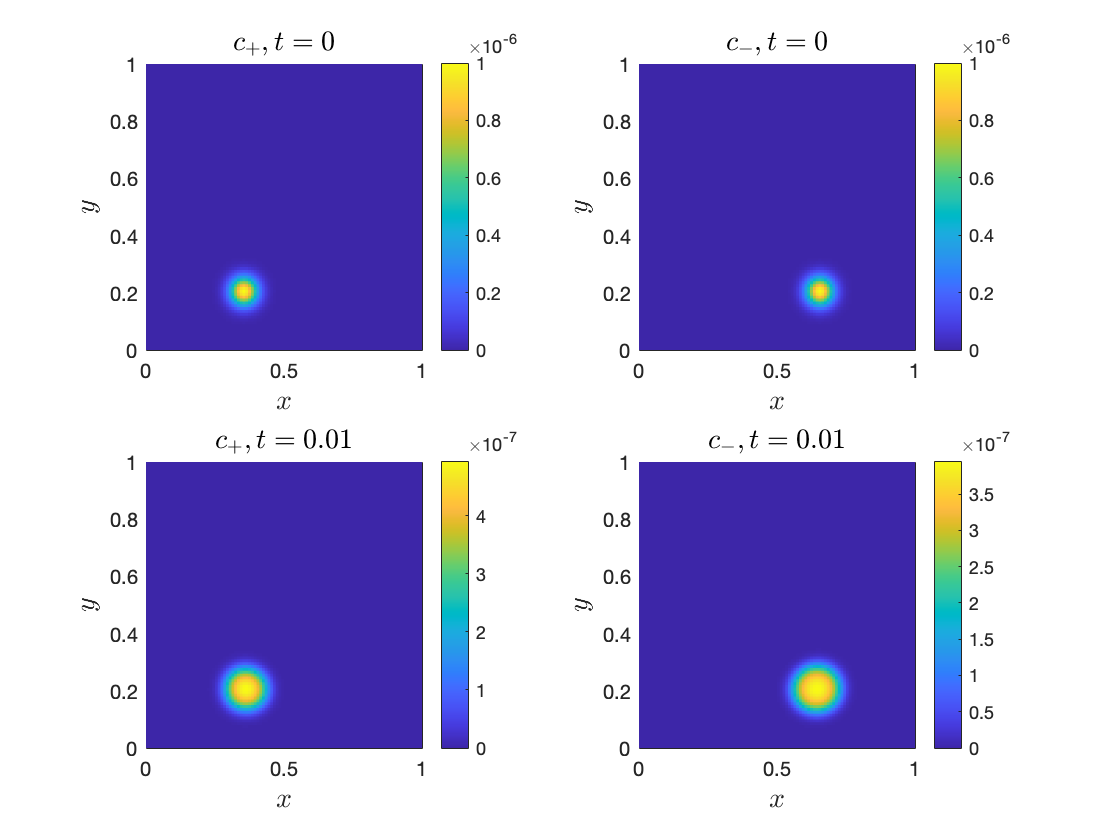}  
\end{overpic}
    \end{minipage}
\begin{minipage}{.49\textwidth}
\begin{overpic}[abs,width=\textwidth,unit=1mm,scale=.25]{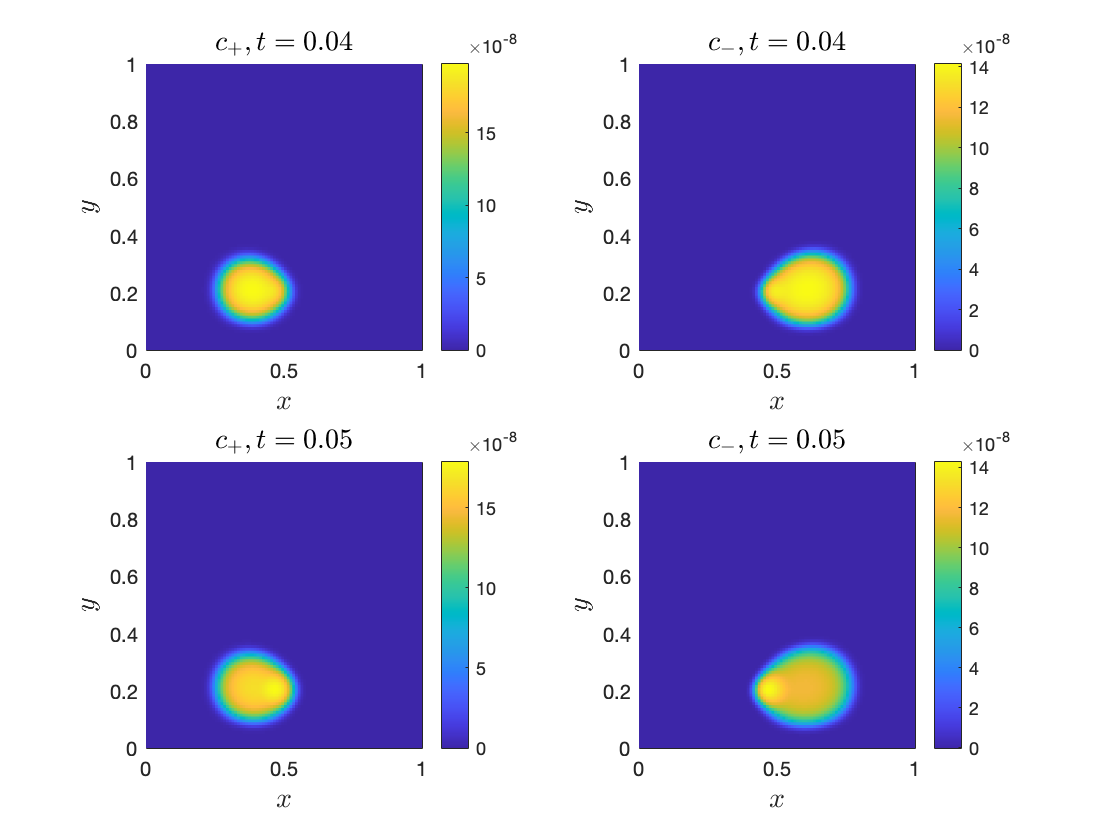}  
\end{overpic}
    \end{minipage}
\begin{minipage}{.49\textwidth}
\begin{overpic}[abs,width=\textwidth,unit=1mm,scale=.25]{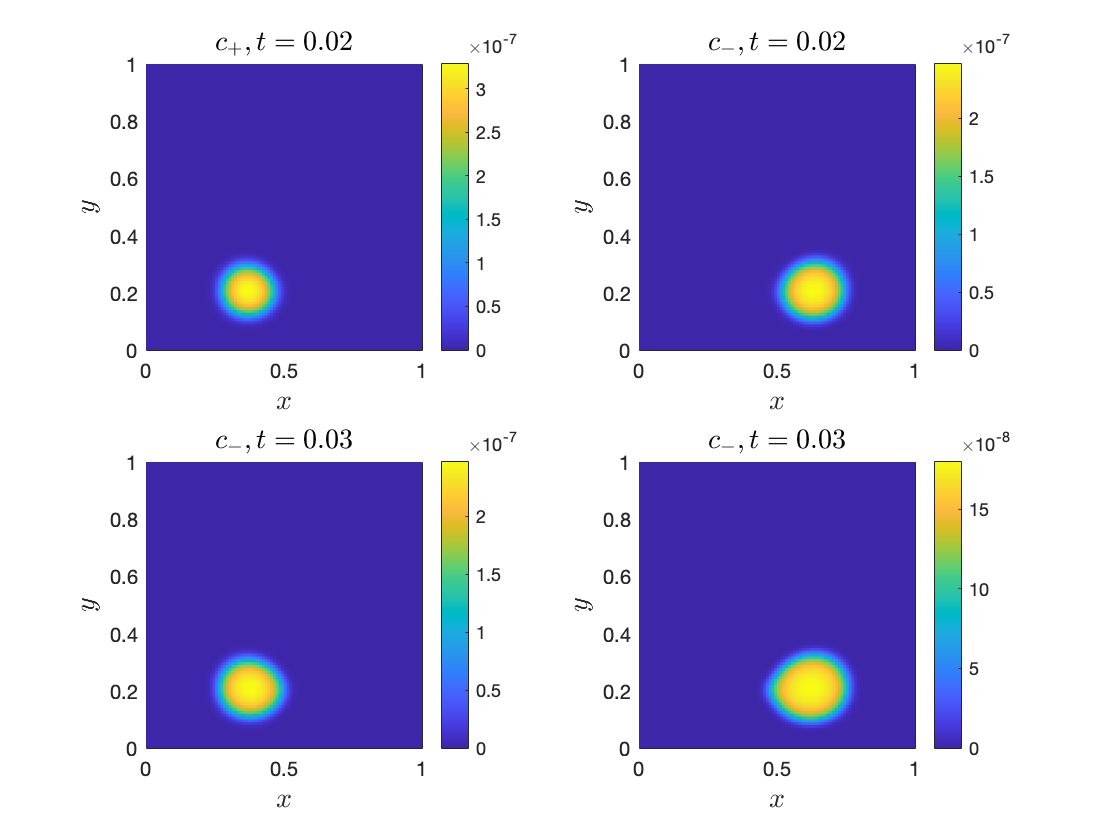}  
\end{overpic}
    \end{minipage}
\begin{minipage}{.49\textwidth}
\begin{overpic}[abs,width=\textwidth,unit=1mm,scale=.25]{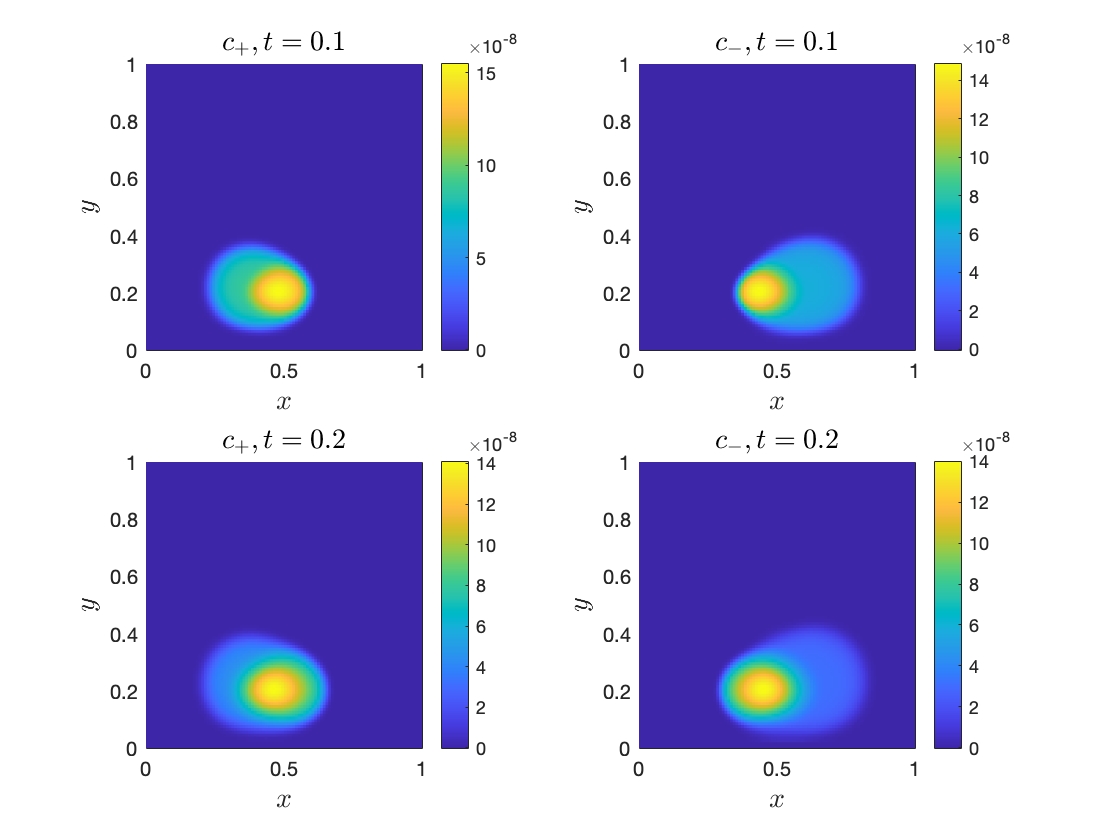}  
\end{overpic}
    \end{minipage}    
\caption{\textit{Time evolution of the system~\eqref{eq:vectorial},\eqref{eq:vectorial2} reaching the quasi-neutrality limit, with $\varepsilon = 10^{-11}$. Time discretization chosen is I2: IMEX-SA(2,2,2), $x^{\rm in}_\pm = 0.5 \mp 0.15$ and $\Delta t = 0.1h$.}}  
\label{fig:solutimes11_2}
\end{figure}

\begin{figure}[h]
\centering
\includegraphics[width=0.6\textwidth]{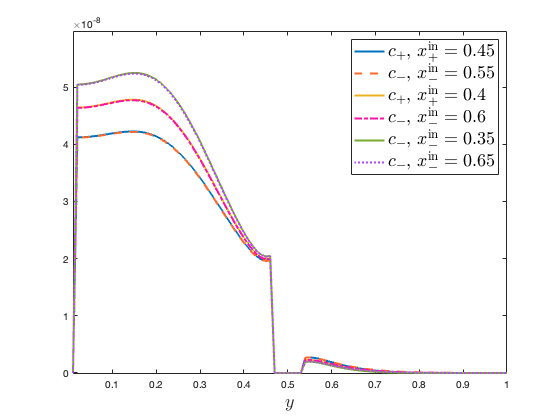}
\caption{\textit{Profiles of the concentrations $c_+$ (solid lines) and $c_-$ (dashed lines) in two dimensions, at $x = 0.5$ and time $t = 10$. We show the solutions for different choices of initial conditions. }}  
\label{fig:profiles}
\end{figure}

Figs.~\ref{fig:solutimes9}-\ref{fig:solutimes11_2} illustrate the time evolution of the system~\eqref{eq:vectorial},\eqref{eq:vectorial2} where the concentrations approach the quasi-neutrality just before starting to diffuse together. Fig.~\ref{fig:solutimes9} corresponds to $\varepsilon = 10^{-9}$, Fig.~\ref{fig:solutimes10} to $\varepsilon = 10^{-10}$, Figs.~\ref{fig:solutimes11_3}-\ref{fig:solutimes11_2} to $\varepsilon = 10^{-11}$, showing how rapidly the solutions reach to quasi-neutrality due to the influence of the Debye length. In Figs.~\ref{fig:solutimes11_3}-\ref{fig:solutimes11_2} we choose different initial conditions to prove the robustness of the numerical scheme in handling different initial data and to show the different transient dynamics in these cases: $x_\pm^{\rm in} = 0.5 \mp 0.05 $ in Fig.~\ref{fig:solutimes11_3}, $x_\pm^{\rm in} = 0.5 \mp 0.1 $ in Fig.~\ref{fig:solutimes11} and $x_\pm^{\rm in} = 0.5 \mp 0.15 $ in Fig.~\ref{fig:solutimes11_2}. As a consequence of considering different initial conditions, Fig.~\ref{fig:profiles} shows the profiles of $c_\pm$ at time $t = 10$. We observe that the larger the difference in the initial peaks of the Gaussians, i.e., $|x^{\mathrm{in}}_+ - x^{\mathrm{in}}_-|$, the later the concentrations begin to merge and to diffuse. This behavior indicates that an increased initial distance between the two species delays their mutual interaction and leads to a slower diffusive dynamics over time.

\begin{table}[h!]
\label{table_parameters}
\centering
\caption{\textit{Parameters involved}} 
\begin{tabular}{cccccc}
\toprule
Symbol & value & Symbol & value & Symbol & value\\ 
\midrule
$D_0$ & $10^{-9}m^2 s^{-1}$ & $D^+/D_0$ & $1.5$  & $D^-/D_0$ & $0.5$\\
$m_0$ & $10^{-3}\,$ Kg$\, mol^{-1}$ & $m^+$ & $23$ & $m^-$ & $265$\\
\bottomrule
\end{tabular}
\end{table}

\section{Conclusions}
In this work, we compare different formulations and time discretizations in order to identify the most suitable approach with respect to asymptotic properties and computational efficiency, depending on the chosen Debye length. While the formulation presented in system~\eqref{eq_full_model_adim} is computationally faster, thanks to the independence of the two parabolic equations and the simplicity of the nonlinear drift term (a single addend), it faces a limitation as the Debye length tends to zero: the severe restriction on the time step makes the numerical scheme inefficient and impractical. These observations suggest that the new formulation in system~\eqref{eq_CQP_system}, being asymptotic preserving, offers a more robust performance for small Debye lengths, ensuring both efficiency and asymptotic accuracy. Furthermore, the numerical scheme shows adaptability to different choices of initial conditions. This property has been verified in our numerical experiments, confirming that the method is stable with different initial configurations.

As a natural continuation of this work, we plan to couple the QNL formulation of the Poisson-Nernst-Planck model with oscillating bubble dynamics, described by the Navier-Stokes equations for the velocity field of the surrounding fluid, in order to investigate ion transport under time-dependent boundary conditions. This extension will allow us to explore the interplay between interfacial motion, surface charge evolution, and the surrounding flow field in more realistic physical configurations. 

\begin{table}[h]
\label{tab:time}
\centering
\caption{\textit{Execution time per iteration for the two systems as $\varepsilon$ varies. With $c_\pm$-system we refer to \eqref{eq:vectorial1}, while with $\mathcal{C},\mathcal{Q}$-system to \eqref{eq:vectorial2}. Time discretization chosen is I2: IMEX-SA(2,2,2).}} 
\begin{tabular}{ccc}
\toprule
$\boldsymbol{\varepsilon}$ & \textbf{$c_\pm$-system} & \textbf{$\mathcal{C},\mathcal{Q}$-system} \\ 
\midrule
$10^{-4}$ & 0.21\,s & 0.23\,s \\ 
$10^{-9}$ & 0.31\,s & 1.88\,s \\ 
\bottomrule
\end{tabular}
\end{table}

\section*{Acknowledgments}
{The work has been supported by the Spoke 10 Future AI Research (FAIR) of the Italian Research Center funded by the Ministry of University and Research as part of the National Recovery and Resilience Plan (PNRR). \\ 
The authors have also been supported also by Italian Ministerial grant PRIN 2022 'Efficient numerical schemes and optimal control methods for time-dependent partial differential equations', No. 2022N9BM3N - Finanziato dall'Unione europea - Next Generation EU.} \\
{The work has also been supported by the Italian Ministerial grant PRIN 2022 PNRR 'FIN4GEO: Forward and Inverse Numerical Modeling of hydrothermal systems in volcanic regions with application to geothermal energy exploitation', No. P2022BNB97 - Finanziato dall'Unione europea - Next Generation EU. } \\
The authors have also been supported also by Italian Ministerial grant PRIN 2022 'Advanced numerical methods for time dependent parametric partial differential equations with applications', No. 2022KA3JBA - Finanziato dall'Unione europea - Next Generation EU. \\
{The authors are members of the Gruppo Nazionale Calcolo Scientifico-Istituto Nazionale di Alta Matematica (GNCS-INdAM).}
We acknowledge the CINECA award under the ISCRA initiative, for the availability of high-performance computing resources.

\section*{Appendix}
\appendix
\subsection*{Details of the space discretization}
\label{sec:app_space_discr}
In this section, we describe some detail of the space discretization employed in the numerical tests. As mentioned in Section~\ref{sec:space_discretization}, we apply a numerical scheme recently proposed in \cite{astuto2025nodal} for elliptic partial differential  equations in arbitrary domains.

The domain $\Omega$ is implicitly defined by the level set function $\phi(x,y)$ that is negative inside $\Omega$, positive in $R\setminus \Omega$ and zero on the boundary $\Gamma$.
In this paper, we choose the signed distance function between $(x,y)$ and $\Gamma$, i.e.,\ $\phi(x,y)=\sqrt{(x-x_c)^2+(y-y_c)^2} - r_\mathcal{B}$, where $r_\mathcal{B}$ is the radius of the circle.

The set of grid points will be denoted by $\mathcal N$, with $\# \mathcal N = (1+N)^2$, the active nodes (i.e.,\ internal $\mathcal{I}$ or ghost $\mathcal{G}$) by $\mathcal A = \mathcal{I}\cup\mathcal{G} \subset \mathcal N$, the set of inactive points by $\mathcal O \subset \mathcal N$, with $\mathcal O\cup\mathcal A = \mathcal N$ and $\mathcal O \cap \mathcal A = \emptyset$ and the set of cells by $\mathcal C$, with $\# \mathcal C = N^2$. Finally, we denote by $\Omega_c = R\setminus \Omega$ the outer region in $R$.

Here, we define the set of ghost points $\mathcal{G}$, which are grid points that belong to $\Omega_c$, and are vertices of cut cells, formally defined as
\begin{equation}
\notag
	(x,y) \in \mathcal{G} \iff (x,y) \in {\mathcal N}\cap \Omega_c  \text{ and } \{(x \pm h,y),(x,y\pm h), (x \pm h,y\pm h) \} \cap \mathcal I \neq \emptyset.
\end{equation}
The discrete spaces $V_h$ and $Q_h$ are given by the piecewise bilinear functions which are continuous in $R$.
As a basis of $V_h$ and $Q_h$, we choose functions whose restriction to the cell is a polynomial in $\mathbb{Q}_1$.


Making use of the basis functions defined in Eq.~\eqref{eq:basis} as test functions in the variational problems, we rewrite (\ref{pro:variational1}-\ref{pro:variational}), as  
\begin{pro}
Find $c_{\pm,h} \in V_h$ and $\Phi_h \in Q_h$ such that,  for almost every $t\in(0,T)$, it holds
\begin{subequations}
\label{pro:variational1_h}
\begin{align}
\left( \frac{\partial c_{+,h}}{\partial t}, v_h\right)_{L^2(\Omega_h)} = & - D_+ \left(\left( \nabla c_{+,h},\nabla v_h \right)_{L^2(\Omega_h)} + \left(c_{+,h} \nabla \Phi_h,\nabla v_h \right)_{L^2(\Omega_h)} \right) \\ 
\left( \frac{\partial c_{-,h}}{\partial t}, u_h\right)_{L^2(\Omega_h)} = & - D_- \left(\left( \nabla c_{-,h},\nabla u_h \right)_{L^2(\Omega_h)} - \left(c_{-,h} \nabla \Phi_h,\nabla u_h \right)_{L^2(\Omega_h)} \right) \\
\varepsilon \left( \nabla \Phi_h,\nabla q_h\right)_{L^2(\Omega_h)} = &\left(\frac{c_{+,h}}{m_+} - \frac{c_{-,h}}{m_-},q_h\right)_{L^2(\Omega_h)}
\end{align}
\end{subequations}
\end{pro}
\begin{figure}[h]
    \centering
    \begin{minipage}{.4\textwidth}
\begin{overpic}[abs,width=0.75\textwidth,unit=1mm,scale=.25]{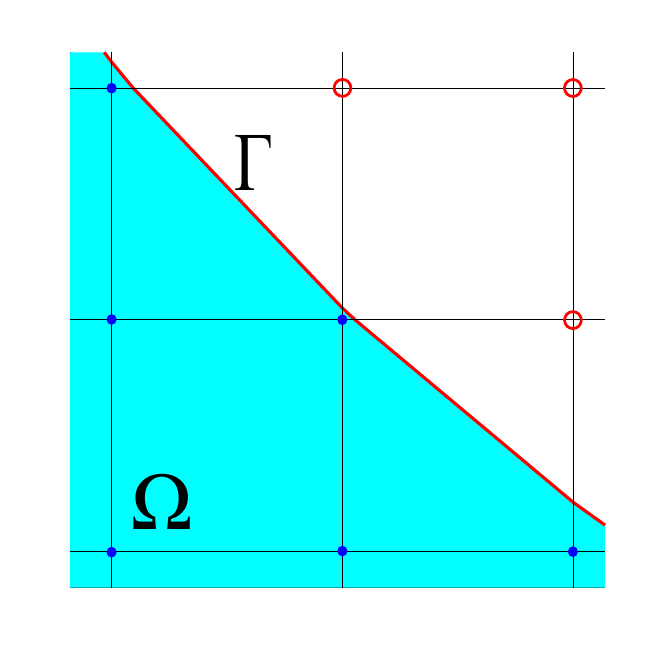}  
\put(-1,33){(a)}
\put(19.5,19.5){P}
\put(14,24.5){$\mathcal B$}
\end{overpic}
    \end{minipage}
    \begin{minipage}{.4\textwidth}
\begin{overpic}[abs,width=0.75\textwidth,unit=1mm,scale=.25]{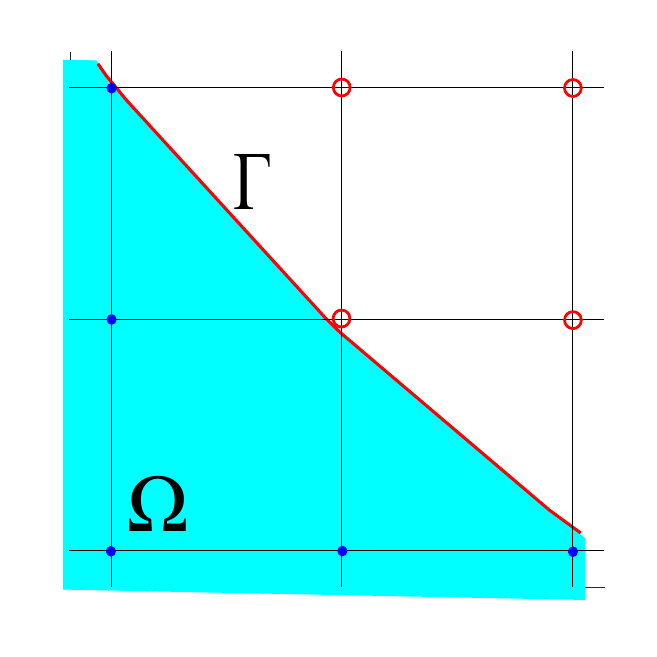}        
\put(-1,33){(b)}
\put(14,23.5){$\mathcal B$}
\put(19.5,19.5){G}
\end{overpic}
    \end{minipage}
\caption{\textit{Grid before and after snapping technique. (a): representation of the cell related to the internal point $P$ (blue points), whose distance from $\Gamma_\mathcal{B}$ is less than $h^2$, i.e. $\phi(P)<h^2$; (b): zoom-in of the shape of the domain, after the grid point $P$ has changed its classification, from internal to ghost point (red circles).
}}  
\label{fig:snapping}
\end{figure}
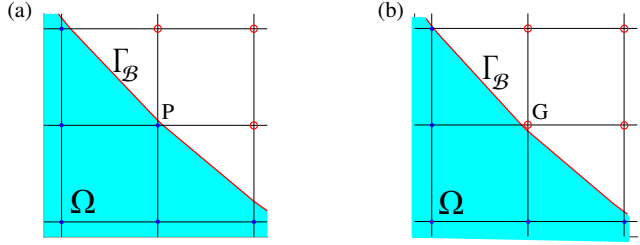

\begin{pro}\label{pro:variational_2}
Find $\mathcal C_h, \mathcal Q_h \in V_h$ and $\Phi_h \in Q_h$  such that, for almost every $t\in(0,T)$, it holds
\begin{subequations}
\label{pro:variational_h}
\begin{align} 
\left( \frac{\partial \mathcal{C}_h}{\partial t}, v_h \right)_{L^2(\Omega_h)} = & - \widetilde D \left( \nabla \mathcal{C}_h,\nabla v_h  \right)_{L^2(\Omega_h)}  - \varepsilon \widehat D \left( \nabla \mathcal Q_h, \nabla v_h  \right)_{L^2(\Omega_h)} \\ \nonumber & - \left(\left( \widehat D \mathcal C_h + \varepsilon \widetilde D \mathcal Q_h \right) \nabla \Phi_h,\nabla v_h  \right)_{L^2(\Omega_h)}  \\
{\varepsilon}\left( \frac{\partial \mathcal{Q}_h}{\partial t}, u_h \right)_{L^2(\Omega_h)} = & - {\widehat D} \left( \nabla \mathcal{C}_h,\nabla u_h  \right)_{L^2(\Omega_h)}  -  {\varepsilon}\widetilde D \left( \nabla \mathcal Q_h, \nabla u_h  \right)_{L^2(\Omega_h)} \\\nonumber  &- \left(\left( {\widetilde D} \mathcal C_h +  {\varepsilon}\widehat D \mathcal Q_h \right) \nabla \Phi_h,\nabla u_h  \right)_{L^2(\Omega_h)}  \\ 
\left( \nabla \Phi_h,\nabla q_h\right)_{L^2(\Omega_h)} = &\left(\mathcal{Q}_h,q_h\right)_{L^2(\Omega_h)}.
\end{align}
\end{subequations}
\end{pro}

When dealing with arbitrary domains, cut cells may arise at the boundary. Since the size of these cells cannot be controlled, they may lead to stability issues by increasing the condition number of the computational matrices. For this reason, particular attention should be devoted to the regions close to the boundary. In \cite{astuto2025nodal}, the authors propose a snapping-back-to-grid strategy that prevents the cut cell size from becoming 'too small', namely smaller than the order of the numerical scheme (i.e., $h^2$). In Fig.~\ref{fig:snapping} (a), illustrates a case where a cut cell is identified as 'too small' (meaning that the value of the level set function in $P$ is less than $h^2$). In such case, the internal point P changes its classification to ghost point G.

\end{document}